\documentclass[10pt]{amsart}
\usepackage{amsmath}
\usepackage{amsxtra}
\usepackage{amscd}
\usepackage{amsthm}
\usepackage{amsfonts}
\usepackage{amssymb}
\usepackage{eucal}

\usepackage{latexsym}

\newtheorem{cor}[subsection]{Corollary}
\newtheorem{lem}[subsection]{Lemma}
\newtheorem{prop}[subsection]{Proposition}

\newtheorem{conj}[subsection]{Conjecture}

\newtheorem{thm}[subsection]{Theorem}

\theoremstyle{remark}

\theoremstyle{definition}

\newcommand{\thmref}[1]{Theorem~\ref{#1}}

\newcommand{\secref}[1]{Sect.~\ref{#1}}
\newcommand{\lemref}[1]{Lemma~\ref{#1}}
\newcommand{\propref}[1]{Proposition~\ref{#1}}
\newcommand{\corref}[1]{Corollary~\ref{#1}}
\newcommand{\conjref}[1]{Conjecture~\ref{#1}}

\newcommand{\nc}{\newcommand}
\nc{\renc}{\renewcommand}
\nc{\ssec}{\subsection}
\nc{\sssec}{\subsubsection}
\nc{\on}{\operatorname}

\nc\ol{\overline}
\nc\ul{\underline}
\nc\wt{\widetilde}
\nc\tboxtimes{\wt{\boxtimes}}
\nc{\wh}{\widehat}
\nc{\mc}{\mathcal}

\nc{\CM}{{\mathcal M}}
\nc{\CN}{{\mathcal N}}
\nc{\CF}{{\mathcal F}}
\nc{\D}{{\mathcal D}}
\nc{\CQ}{{\mathcal Q}}
\nc{\CY}{{\mathcal Y}}
\nc{\CX}{{\mathcal X}}
\nc{\CG}{{\mathcal G}}
\nc{\CE}{{\mathcal E}}
\nc{\CC}{{\mathcal C}}
\nc{\CO}{{\mathcal O}}
\renc{\CC}{{\mathcal C}}
\nc{\CT}{{\mathcal T}}
\nc{\CK}{{\mathcal K}}
\nc{\CS}{{\mathcal S}}
\nc{\CH}{{\mathcal H}}
\nc{\CU}{{\mathcal U}}
\nc{\CV}{{\mathcal V}}
\nc{\CA}{{\mathcal A}}
\nc{\CB}{{\mathcal B}}
\nc{\CW}{{\mathcal W}}
\nc{\CL}{{\mathcal L}}
\nc{\CP}{{\mathcal P}}
\nc{\CI}{{\mathcal I}}
\nc{\CJ}{{\mathcal J}}
\nc{\CR}{{\mathcal R}}

\nc{\BA}{{\mathbb{A}}}
\nc{\BC}{{\mathbb{C}}}
\nc{\BG}{{\mathbb{G}}}
\nc{\BM}{{\mathbb{M}}}
\nc{\BN}{{\mathbb{N}}}
\nc{\BP}{{\mathbb{P}}}
\nc{\BR}{{\mathbb{R}}}
\nc{\BZ}{{\mathbb{Z}}}
\nc{\BV}{{\mathbb{V}}}
\nc{\BW}{{\mathbb{W}}}
\nc{\BS}{{\mathbb{S}}}
\nc{\BD}{{\mathbb{D}}}
\nc{\BQ}{{\mathbb{Q}}}
\nc{\BL}{{\mathbb{L}}}
\renc{\BW}{{\mathbb{W}}}

\nc{\fa}{{\mathfrak{a}}}
\nc{\fb}{{\mathfrak{b}}}
\nc{\fg}{{\mathfrak{g}}}
\nc{\fgl}{{\mathfrak{gl}}}
\nc{\fh}{{\mathfrak{h}}}
\nc{\fj}{{\mathfrak{j}}}
\nc{\fm}{{\mathfrak{m}}}
\nc{\fl}{{\mathfrak{l}}}
\nc{\fn}{{\mathfrak{n}}}
\nc{\fu}{{\mathfrak{u}}}
\nc{\fp}{{\mathfrak{p}}}
\nc{\fr}{{\mathfrak{r}}}
\nc{\fs}{{\mathfrak{s}}}
\nc{\fsl}{{\mathfrak{sl}}}

\nc{\hsl}{{\widehat{\mathfrak{sl}}}}
\nc{\hgl}{{\widehat{\mathfrak{gl}}}}
\nc{\hg}{{\widehat{\mathfrak{g}}}}
\nc{\hb}{{\widehat{\mathfrak{b}}}}
\nc{\hn}{{\widehat{\mathfrak{n}}}}

\nc{\fA}{{\mathfrak{A}}}
\nc{\fB}{{\mathfrak{B}}}
\nc{\fO}{{\mathfrak{O}}}
\nc{\fD}{{\mathfrak{D}}}
\nc{\fE}{{\mathfrak{E}}}
\nc{\fF}{{\mathfrak{F}}}
\nc{\fG}{{\mathfrak{G}}}
\nc{\fK}{{\mathfrak{K}}}
\nc{\fL}{{\mathfrak{L}}}
\nc{\fC}{{\mathfrak{C}}}
\nc{\fM}{{\mathfrak{M}}}
\nc{\fN}{{\mathfrak{N}}}
\nc{\fH}{{\mathfrak{H}}}
\nc{\fP}{{\mathfrak{P}}}
\nc{\fU}{{\mathfrak{U}}}
\nc{\fV}{{\mathfrak{V}}}
\nc{\fZ}{{\mathfrak{Z}}}
\nc{\fz}{{\mathfrak{z}}}

\nc{\bc}{{\mathbf{c}}}
\nc{\bd}{{\mathbf{d}}}
\nc{\bh}{{\mathbf{h}}}
\nc{\be}{{\mathbf{e}}}
\nc{\bj}{{\mathbf{j}}}
\nc{\bn}{{\mathbf{n}}}
\nc{\bp}{{\mathbf{p}}}
\nc{\bg}{{\mathbf{g}}}
\nc{\bq}{{\mathbf{q}}}
\nc{\bs}{{\mathbf{s}}}
\nc{\bu}{{\mathbf{u}}}
\nc{\bv}{{\mathbf{v}}}
\nc{\bx}{{\mathbf{x}}}
\nc{\by}{{\mathbf{y}}}
\nc{\bw}{{\mathbf{w}}}
\nc{\bA}{{\mathbf{A}}}
\nc{\bK}{{\mathbf{K}}}
\nc{\bB}{{\mathbf{B}}}
\nc{\bC}{{\mathbf{C}}}
\nc{\bD}{{\mathbf{D}}}
\nc{\bH}{{\mathbf{H}}}
\nc{\bM}{{\mathbf{M}}}
\nc{\bN}{{\mathbf{N}}}
\nc{\bV}{{\mathbf{V}}}
\nc{\bW}{{\mathbf{W}}}
\nc{\bL}{{\mathbf{L}}}
\nc{\bU}{{\mathbf{U}}}
\nc{\bX}{{\mathbf{X}}}
\nc{\bI}{{\mathbf{I}}}
\nc{\bZ}{{\mathbf{Z}}}
\nc{\bS}{{\mathbf{S}}}

\nc{\sA}{{\mathsf{A}}}
\nc{\sB}{{\mathsf{B}}}
\nc{\sC}{{\mathsf{C}}}
\nc{\sD}{{\mathsf{D}}}
\nc{\sF}{{\mathsf{F}}}
\nc{\sH}{{\mathsf{H}}}
\nc{\sG}{{\mathsf{G}}}
\nc{\sK}{{\mathsf{K}}}
\nc{\sM}{{\mathsf{M}}}
\nc{\sO}{{\mathsf{O}}}
\nc{\sQ}{{\mathsf{Q}}}
\nc{\sP}{{\mathsf{P}}}
\nc{\sV}{{\mathsf{V}}}
\nc{\sZ}{{\mathsf{Z}}}
\nc{\sfp}{{\mathsf{p}}}
\nc{\sr}{{\mathsf{r}}}
\nc{\sg}{{\mathsf{g}}}
\nc{\sk}{{\mathsf{k}}}
\nc{\ssf}{{\mathsf{f}}}
\nc{\ssh}{{\mathsf{h}}}
\nc{\sse}{{\mathsf{e}}}
\nc{\sfb}{{\mathsf{b}}}
\nc{\sfc}{{\mathsf{c}}}
\nc{\sd}{{\mathsf{d}}}

\nc{\Av}{\on{Av}}
\nc{\act}{\on{act}}
\nc{\Hom}{\on{Hom}}
\nc{\End}{\on{End}}
\nc{\Lie}{\on{Lie}}
\nc{\Loc}{\on{Loc}}
\nc{\IC}{\on{IC}}
\nc{\Aut}{\on{Aut}}
\nc{\rk}{\on{rk}}
\nc{\Sh}{\on{Sh}}
\nc{\Perv}{\on{Perv}}
\nc{\pos}{{\on{pos}}}
\nc{\Conv}{\on{Conv}}
\nc{\Sph}{\on{Sph}}
\nc{\Sym}{\on{Sym}}
\nc{\Rep}{\on{Rep}}
\nc{\RepH}{{\mc R}ep(H)}
\nc{\Fun}{\on{Fun}}
\nc{\Id}{\on{Id}}
\nc{\id}{\on{id}}
\renc{\mod}{\on{--mod}}

\nc{\oG}{\overset{\circ}{G}{}}
\nc{\oGB}{{\overset{\circ}{G/B}{}}}
\nc{\oGN}{{\overset{\circ}{G/N}{}}}
\nc{\uBC}{\underline{\BC}}

\nc{\crit}{{\on{crit}}}
\nc{\reg}{{\on{reg}}}
\nc{\nilp}{{\on{nilp}}}
\nc{\ord}{\on{ord}}
\nc{\nil}{\wt{\on{reg}}}
\nc{\mb}{\mathbf}
\nc{\ren}{\on{ren}}
\nc{\res}{\on{res}}
\nc{\RS}{{\on{RS}}}
\nc{\Dist}{\on{Dist}}
\nc{\semiinf}{{\frac{\infty}{2}}}
\nc{\semiinfi}{{\frac{\infty}{2}+i}}
\nc{\semiinfb}{{\frac{\infty}{2}+\bullet}}
\nc{\torsemiinf}{{\overset{\semiinf}\otimes}}
\nc{\Hitch}{\on{Hitch}}

\nc{\hl}{\overset{\leftarrow}h}
\nc{\hr}{\overset{\rightarrow}h}
\nc\Dh{\widehat{\D}}
\nc{\Gr}{\on{Gr}}
\nc{\Grb}{\ol{\Gr}{}}
\nc{\Fl}{\on{Fl}}
\nc{\Flt}{\wt{\Fl}{}}
\nc{\Pic}{\on{Pic}}
\nc{\Bun}{\on{Bun}}

\nc{\bDR}{\mathbf {DR}}
\nc{\uV}{\underline{V}}
\nc{\arrowtimes}{\overset{\to}\otimes}
\nc{\hattimes}{\widehat\otimes}
\nc{\larrowtimes}{\overset{\leftarrow}\otimes}
\nc{\shriektimes}{\overset{!}\otimes}
\nc{\startimes}{\overset{*}\otimes}
\nc{\sCliff}{\mathsf {Cliff}}
\nc{\sSpin}{\mathsf {Spin}}

\nc{\one}{{\mathbf{1}}}

\nc\Spec{\on{Spec}}
\nc{\Pro}{\on{Pro}}
\nc{\QCoh}{\on{QCoh}}
\nc{\uHom}{\underline{\on{Hom}}}
\nc{\RHom}{\on{RHom}}
\nc{\uRHom}{\underline{\on{RHom}}}
\nc{\CHom}{{\mathcal Hom}}
\nc{\uCHom}{\underline{{\mathcal Hom}}}
\nc{\uCRHom}{\underline{{\mathcal R}{\mathcal Hom}}}

\nc{\cg}{{\check \fg}}
\nc{\Op}{\on{Op}}
\nc{\nOp}{\on{Op}^{\nilp}_{\cg}}
\nc{\nMOp}{\on{MOp}^{\nilp}_{\cg}}
\nc{\rOp}{\on{Op}^{\reg}_{\cg}}

\nc{\tg}{\wt{\check \fg}}
\nc{\cn}{\check \fn}
\nc{\tn}{\wt{\cn}}
\nc{\cG}{{\check G}}
\nc{\cB}{\check B}
\nc{\cT}{\check T}
\nc{\cH}{\check H}
\nc{\cb}{\check \fb}
\nc{\cN}{\check N}
\nc{\MOp}{\on{MOp}}
\nc{\tN}{\wt{\CN}_{\cG}}
\nc{\dIsom}{{\mathsf{Isom}}_{\Op}}
\nc{\disom}{{\mathsf{isom}}_{\Op}}
\nc{\Kdv}{{\mathsf{Isom}}_{\Op^\reg}}
\nc{\kdv}{{\mathsf{isom}}_{\Op^\reg}}
\nc{\Isom}{{\mathsf{Isom}}}
\nc{\isom}{{\mathsf{isom}}}

\nc{\wcosta}{j_{\wt{w},*}}
\nc{\wsta}{j_{\wt{w},!}}
\nc{\wcost}{j_{w,*}}
\nc{\wst}{j_{w,!}}

\nc{\epsi}{{\mathbf e}^\psi}
\nc{\epsip}{{\mathbf e}^{\psi'}}

\nc{\Ppi}{{\mathbf \Pi}}

\nc{\hCO}{{\hat{\CO}}}
\nc{\hCK}{{\hat{\CK}}}

\nc{\CPreg}{\CP_{G,\on{Op}^\reg}}
\nc{\CPBreg}{\CP_{B,\on{Op}^\reg}}
\nc{\CPnilp}{\CP_{G,\on{Op}^\nilp}}
\nc{\CPBnilp}{\CP_{B,\on{Op}^\nilp}}
\nc{\CPla}{\CP_{G,\on{Op}_{\cla}}}
\nc{\CPBla}{\CP_{B,\on{Op}_{\cla}}}

\nc{\Cat}{\hg_\crit\mod^{I,m}_\nilp}
\nc{\Catf}{{}^f\hg_\crit\mod^{I,m}_\nilp}
\nc{\DCat}{D^b(\hg_\crit\mod_\nilp)^{I^0}}
\nc{\DCatf}{{}^f D^b(\hg_\crit\mod_\nilp)^{I^0}}
\nc{\Catr}{\hg_\crit\mod^{I,m}_\reg}
\nc{\Catrf}{{}^f\hg_\crit\mod^{I,m}_\reg}
\nc{\DCatr}{D^b(\hg_\crit\mod_\reg)^{I^0}}
\nc{\DCatrf}{{}^f D^b(\hg_\crit\mod_\reg)^{I^0}}

\nc{\ch}{\mbox{ch}}
\nc{\Z}{{\mathbb Z}}
\nc{\C}{{\mathbb C}}
\nc{\pone}{{\mathbb C}{\mathbb P}^1}
\nc{\pa}{\partial}
\nc{\F}{{\mathcal F}}
\nc{\arr}{\rightarrow}
\nc{\larr}{\longrightarrow}
\nc{\al}{\alpha}
\nc{\ri}{\rangle}
\nc{\lef}{\langle}
\nc{\W}{{\mathcal W}}
\nc{\la}{\lambda}
\nc{\ep}{\epsilon}
\nc{\su}{\widehat{{\mathfrak s}{\mathfrak l}}_2}
\nc{\sw}{{\mathfrak s}{\mathfrak l}}
\nc{\g}{{\mathfrak g}}
\nc{\h}{{\mathfrak h}}
\nc{\n}{{\mathfrak n}}
\nc{\N}{\widehat{\n}}
\nc{\De}{\Delta}
\nc{\gt}{\widetilde{\g}}
\nc{\Ga}{\Gamma}
\nc{\z}{{\mathfrak Z}}
\nc{\La}{\Lambda}
\nc{\cri}{_{\kappa_c}}
\nc{\kk}{h^\vee}
\nc{\sun}{\widehat{\sw}_N}
\nc{\si}{\sigma}
\nc{\el}{\ell}
\nc{\bi}{\bibitem}
\nc{\om}{\omega}
\nc{\ds}{\displaystyle}
\nc{\dzz}{\frac{dz}{z}}
\nc{\Res}{\on{Res}}
\nc{\Cal}{\mathcal}
\nc{\bb}{{\mathfrak b}}
\nc{\ot}{\otimes}
\nc{\R}{{\mc R}}
\nc{\yy}{{\mc Y}}
\nc{\ga}{\gamma}

\nc{\us}{\underset}
\nc{\opl}{\oplus}
\nc{\beq}{\begin{equation}}
\nc{\Fq}{{\mathbb F}_q}
\nc{\Mq}{{\mathcal M}}
\nc{\lan}{\langle}
\nc{\ran}{\rangle}

\nc{\Vect}{\on{Vect}}
\nc{\ghat}{\wh\fg}
\nc{\T}{\mc T}
\nc{\Tloc}{\T^\g_{\on{loc}}}
\nc{\vac}{|0\ran}
\nc{\Wick}{{\mb :}}
\nc{\delz}{\partial_z}
\nc{\K}{{\cali K}}
\nc{\cali}{\mathcal}
\nc{\li}{\mathfrak l}
\nc{\lt}{\widetilde{\li}}
\nc{\astar}{a^*}
\nc{\cA}{{\mc A}}
\nc{\ka}{\kappa}

\nc{\OO}{{\mc O}}
\nc{\AutO}{\on{Aut}\OO}
\nc{\DerO}{\on{Der}\OO}
\nc{\DerpO}{\on{Der}_+\OO}
\nc{\Au}{{\mc A}ut}
\nc{\mf}{\mathfrak}
\nc{\V}{{\mc V}}
\nc{\hh}{\wh{\h}}

\nc{\pp}{{\mathfrak p}}
\nc{\mm}{{\mathfrak m}}
\nc{\rr}{{\mathfrak r}}
\nc{\ket}{\rangle}
\nc{\zz}{{\mathfrak z}}
\nc{\gr}{\on{gr}}
\nc{\Spe}{\on{Spec}}
\nc{\rv}{\crho}
\nc{\can}{\on{can}}
\nc{\Db}{{\mathbb D}}
\nc{\ww}{w}

\nc{\RR}{\on{R}}
\nc{\PPi}{{\mathbf \Pi}}
\nc{\M}{{\mathbb M}}
\nc{\Mv}{{\mathbb M}^\vee}
\nc{\VV}{{\mathbb V}}
\nc{\bsl}{\backslash}

\nc{\bchi}{{\mathbf {\chi}}}
\nc{\anch}{{\mathbf {anch}}}

\nc{\cla}{{\check{\la}}}
\nc{\cmu}{{\check{\mu}}}
\nc{\crho}{{\check{\rho}}}
\nc{\com}{{\check{\omega}}}
\nc{\DD}{{\mc D}}
\nc{\E}{{\mc E}}
\nc{\Ll}{{\mc L}}

\nc{\Conn}{\on{Conn}_{\cH}(\omega^{\crho})}
\nc{\ConnD}{\on{Conn}_H(\omega_{\D}^{\crho})}
\nc{\ConHD}{\on{Conn}_{\check{H}}(\omega_{\D}^{\rho})}
\nc{\ConnDt}{\on{Conn}_H(\omega_{\D^\times}^{\crho})}
\nc{\ConHDt}{\on{Conn}_{\check{H}}(\omega_{\D^\times}^{\rho})}

\nc{\Hecke}{{\on{Hecke}}}

\nc{\cLambda}{{\check\Lambda}}
\nc{\cnu}{{\check\nu}}
\nc{\ceta}{{\check\eta}}

\nc{\Ind}{\on{Ind}}

\nc{\CTop}{{\mathcal Top}}

\nc{\ppart}{(\!(t)\!)}

\nc{\qu}{/\!/}

\nc{\gen}{{gen}}

\nc{\Ext}{\on{Ext}}

\nc{\aff}{{\on{aff}}}

\begin{document}

\title{Localization of $\ghat$--modules on the affine Grassmannian}

\author{Edward Frenkel}\thanks{The research of E.F. was supported by
the DARPA grant HR0011-04-1-0031 and by the NSF grant DMS-0303529.}

\address{Department of Mathematics, University of California,
  Berkeley, CA 94720, USA}

\email{frenkel@math.berkeley.edu}

\author{Dennis Gaitsgory}

\address{Department of Mathematics, Harvard University,
Cambridge, MA 02138, USA}

\email{gaitsgory@math.harvard.edu}

\date{December 2005}

\begin{abstract}

We consider the category of modules over the affine Kac-Moody algebra
$\hg$ of critical level with regular central character. In our
previous paper \cite{FG2} we conjectured that this category is
equivalent to the category of Hecke eigen-D-modules on the affine
Grassmannian $G\ppart/G[[t]]$. This conjecture was motivated by our
proposal for a local geometric Langlands correspondence.  In this
paper we prove this conjecture for the corresponding $I^0$-equivariant
categories, where $I^0$ is the radical of the Iwahori subgroup of
$G\ppart$. Our result may be viewed as an affine analogue of the
equivalence of categories of $\g$-modules and D-modules on the flag
variety $G/B$, due to Beilinson-Bernstein and Brylinski-Kashiwara.

\end{abstract}

\maketitle

\section*{Introduction}

\ssec{} Let $G$ be a simple complex algebraic group and $B$ its Borel
subgroup. Consider the category $\on{D}(G/B)\mod$ of left D-modules on
the flag variety $G/B$. The Lie algebra $\fg$ of $G$, and hence its
universal enveloping algebra $U(\fg)$, acts on the space
$\Gamma(G/B,{\mc F})$ of global sections of any D-module ${\mc
F}$. The center $Z(\fg)$ of $U(\fg)$ acts on $\Gamma(G/B,{\mc F})$ via
the augmentation character $\chi_0: Z(\fg) \to \C$. Let
$\fg\mod_{\chi_0}$ be the category of $\fg$-modules on which $Z(\fg)$
acts via the character $\chi_0$. Thus, we obtain a functor $$\Gamma:
\on{D}(G/B)\mod \to \fg\mod_{\chi_0}.$$ In \cite{BB} A. Beilinson and
J. Bernstein proved that this functor is an equivalence of
categories. Moreover, they generalized this equivalence to the case of
twisted D-modules, for twistings that correspond to dominant weights
$\lambda \in \fh^*$.

\medskip

Let $N$ be the unipotent radical of $B$. We can consider the
$N$--equivariant subcategories on both sides of the above
equivalence. On the D-module side this is the category
$\on{D}(G/B)\mod^N$ of $N$--equivariant D-modules on $G/B$, and on the
$\fg$-module side this is the block of the category ${\mc O}$
corresponding to the central character $\chi_0$. The resulting
equivalence of categories, which follows from \cite{BB}, and which was
proved independently by J.-L. Brylinski and M. Kashiwara \cite{BK}, is
very important in applications to representation theory of $\fg$.

\medskip

Now let $\ghat$ be the affine Kac-Moody algebra, the universal central
extension of the formal loop agebra $\fg\ppart$. Representations of
$\wh\fg$ have a parameter, an invariant bilinear form on $\fg$, which
is called the level. There is a unique inner product $\kappa_{\can}$
which is normalized so that the square length of the maximal root of
$\fg$ is equal to $2$. Any other inner product is equal to
$\kappa = k\cdot \ka_{\can}$, where $k \in \C$, and so a level
corresponds to a complex number $k$. In particular, it makes sense to
speak of {\em integral levels}. Representations, corresponding to the
bilinear form which is equal to minus one half of the Killing form
(for which $k=-h^\vee$, minus the dual Coxeter number of $\g$) are
called representations of {\em critical level}. This is really the
``middle point'' amongst all levels (and not the zero level, as one
might naively expect).

\medskip

There are several analogues of the flag variety in the affine case. In
this paper (except in the Appendix) we will consider exclusively
the {\em affine Grassmannian} $\Gr_G = G\ppart/G[[t]]$.

Another possibility is to consider the affine flag scheme $\Fl_G =
G\ppart/I$, where $I$ is the Iwahori subgroup of $G\ppart$. Most of
the results of this paper that concern the critical level have
conjectural counterparts for the affine flag variety, but they are
more difficult to formulate. In particular, one inevitably has to
consider derived categories, whereas for the affine Grassmannian
abelian categories suffice. We refer the reader to the Introduction of
our previous paper \cite{FG2} for more details.

\medskip

There is a canonical line bundle $\Ll_{\can}$ on $\Gr_G$ such that the
action of $\fg\ppart$ on $\Gr_G$ lifts to an action of
$\hg_{\kappa_{\can}}$ on $\Ll_{\can}$. For each level $\kappa$ we can
consider the category $\on{D}(\Gr_G)_\kappa\mod$ of right D-modules on
$\Gr_G$ twisted by $\Ll_{\can}^{\otimes k}$, where $\kappa=k\cdot
\kappa_{\can}$.  (Recall that although the line bundle $\Ll^{\otimes
k}_{\can}$ only makes sense when $k$ is integral, the corresponding
category of twisted D-module is well-defined for an arbitrary $k$.)
Since $\Gr_G$ is an ind-scheme, the definition of these categories
requires some care (see \cite{BD} and \cite{FG1}).

Let $\hg_\kappa\mod$ be the category of (discrete) modules over
the affine Kac-Moody algebra of level $\kappa$ (see
\secref{recol}). Using the fact that the action of $\fg\ppart$ on
$\Gr_G$ lifts to an action of $\hg_{\kappa_{\can}}$ on $\Ll_{\can}$,
we obtain that for each level $\kappa$ we have a naturally defined
functor of global sections:
\begin{equation} \label{glob sections}
\Gamma: \on{D}(\Gr_G)_\kappa\mod \to \ghat_\kappa\mod.
\end{equation}

The question that we would like to address in this paper is whether
this functor is an equivalence of categories, as in the
finite-dimensional case.

\ssec{}   \label{neg level}

The first results in this direction were obtained in \cite{BD,FG1}.
Namely, in {\it loc. cit.} it was shown that if $\kappa$ is such that
$\kappa=k\cdot \kappa_{\can}$ with $k+h^\vee\notin \BQ^{>0}$, then the
functor $\Gamma$ of \eqref{glob sections} is exact and faithful. (In
contrast, it is known that this functor is not exact for $k+h^\vee\in
\BQ^{>0}$.)  The condition $k+h^\vee\notin \BQ^{>0}$ is analogous to
the dominant weight condition of \cite{BB}.

\medskip

Let us call $\kappa$ {\it negative} if $k+h^\vee\notin \BQ^{\geq
0}$. In this case one can show that the functor of \eqref{glob
sections} is fully faithful. In fact, in this case it makes more sense
to consider $H$-monodromic twisted D-modules on the enhanced affine
flag scheme $\wt\Fl_G = G\ppart/I^0$, rather than simply twisted
D-modules on $\Gr_G$, and the corresponding functor $\Gamma$ to
$\ghat_\kappa\mod$.  The above exactness and fully-faithfulness
assertions are still valid in this context. However, the above functor
is not an equivalence of categories. Namely, the RHS of \eqref{glob
sections} has "many more" objects than the LHS.

When $\kappa$ is integral, A.~Beilinson has proposed a conjectural
intrinsic description of the image of the category
$\on{D}(\wt\Fl_G)_\kappa\mod$ inside $\ghat_\kappa\mod$ (see Remark
(ii) in the Introduction of \cite{Bei}). As far as we know, no such
description was proposed when $\kappa$ is not integral.

\medskip

It is possible, however, to establish a partial result in this
direction.  Namely, let $I^0\subset I$ be the unipotent radical of the
Iwahori subgroup $I$. We can consider the category
$\on{D}(\wt\Fl_G)_\kappa\mod^{I^0}$ of $I^0$-equivariant twisted
D-modules on $\wt\Fl_G$. The corresponding functor $\Gamma$ of global
sections takes values in the affine version of category $\CO$, i.e.,
in the subcategory $\ghat_\kappa\mod^{I^0}\subset \ghat_\kappa\mod$,
whose objects are $\ghat_\kappa$-modules on which the action of the
Lie algebra $\on{Lie}(I^0) \subset \hg_\kappa$ integrates to an
algebraic action of the group $I^0$.

One can show that the functor $\Gamma$ induces an equivalence between
an appropriately defined subcategory of $H$-monodromic objects of
$\on{D}(\wt\Fl_G)_\kappa\mod^{I^0}$ and a specific block of
$\ghat_\kappa\mod^{I^0}$.  This result, which is well-known to
specialists, is not available in the published literature. For the
sake of completeness, we sketch one of the possible proofs in the
Appendix of this paper.

\ssec{}

In this paper we shall concentrate on the case of the critical level,
when $k=-h^\vee$. We will see that this case is dramatically different
from the cases considered above. In \cite{FG2} we made a precise
conjecture describing the relationship between the corresponding
categories $\on{D}(\Gr_G)_\crit\mod$ and $\ghat_\crit\mod$. We shall
now review the statement of this conjecture.

\medskip

First, let us note that the image of the functor $\Gamma$ is in
a certain subcategory of $\ghat_\crit\mod$, singled out by the
condition on the action of the center.

Let $\fZ_\fg$ denote the center of the category $\ghat_\crit\mod$
(which is the same as the center of the completed enveloping algebra
of $\hg_\crit$). The fact that this center is non-trivial is what
distinguishes the critical level from all other levels. Let
$\fZ^\reg_\fg$ denote the quotient of $\fZ_\fg$, through which it acts
on the vacuum module
$\BV_\crit:=\on{Ind}^{\hg_\crit}_{\fg[[t]] \oplus \BC}(\BC)$.

Let $\ghat_\crit\mod_\reg$ be the full subcategory of
$\ghat_\crit\mod$, whose objects are $\ghat_\crit$-modules on which
the action of the center $\fZ_\fg$ factors through $\fZ^\reg_\fg$. It
is known (see \cite{FG1}) that for any $\CF \in
\on{D}(\Gr_G)_\crit\mod$, the space of global sections
$\Gamma(\Gr_G,\CF)$ is an object of $\ghat_\crit\mod_\reg$. (Here and
below we write $M \in {\mc C}$ if $M$ is an object of a category ${\mc
C}$.) Thus, $\ghat_\crit\mod_\reg$ is the category that may be viewed
as an analogue of the category $\fg\mod_{\chi_0}$ appearing on the
representation theory side of the Beilinson-Bernstein equivalence.

\medskip

However, the functor of global sections $\Gamma:
\on{D}(\Gr_G)_\crit\mod \to \ghat_\crit\mod_\reg$ is not full, and
therefore cannot possibly be an equivalence.  The origin of the
non-fullness of $\Gamma$ is two-fold, with one ingredient rather
elementary, and another less so.

First, the category $\ghat_\crit\mod_\reg$ has a large center, namely,
the algebra $\fZ^\reg_\fg$ itself, while the center of the category
$\on{D}(\Gr_G)_\crit\mod$ is the group algebra of the finite group
$\pi_1(G)$ (i.e., it has a basis enumerated by the connected
components of $\Gr_G$).

\medskip

Second, the category $\on{D}(\Gr_G)_\crit\mod$
carries an additional symmetry, namely, an action of the tensor
category $\Rep(\cG)$ of the Langlands dual group $\cG$,
and this action trivializes under the functor $\Gamma$.

In more detail, let us recall that, according to \cite{FF,F:wak}, we
have a canonical isomorphism between $\Spec(\fZ^\reg_\fg)$ and the
space $\Op_\cg(\D)$ of $\cg$-opers on the formal disc $\D$ (we refer
the reader to Sect. 1 of \cite{FG2} for the definition and a detailed
review of opers). By construction, over the scheme $\Op_\cg(\D)$
there exists a canonical principal $\cG$-bundle, denoted by
$\CP_{\cG,\Op}$. Let $\CP_{\cG,\fZ}$ bethe $\cG$-bundle over
$\Spec(\fZ^\reg_\fg)$ corresponding to it under the above
isomorphism. For an object $V \in \Rep(\cG)$ let us denote by
$\CV_\fZ$ the associated vector bundle over $\Spec(\fZ^\reg_\fg)$,
i.e., $\CV_{\fZ} = \CP_{\cG,\fZ} \underset{\cG}\times V$.

Consider now the category $\on{D}(\Gr_G)_\crit\mod^{G[[t]]}$. By
\cite{MV}, this category has a canonical tensor structure, and
as such it is equivalent to the category $\Rep(\cG)$ of algebraic
representations of $\cG$; we shall denote by
$$V\mapsto \CF_V:\Rep(\cG)\to \on{D}(\Gr_G)_\crit\mod^{G[[t]]}$$
the corresponding functor. Moreover, we have a canonical action of
$\on{D}(\Gr_G)_\crit\mod^{G[[t]]}$ as a tensor category on 
$\on{D}(\Gr_G)_\crit\mod$ by convolution functors, $\F \mapsto \F
\star \F_V$.

A. Beilinson and V. Drinfeld \cite{BD} have proved that there
are functorial isomorphisms
$$
\Gamma(\Gr_G,\F \star \F_V) \simeq \Gamma(\Gr_G,\F)
\underset{\fZ^\reg_\fg}\otimes \CV_{\fZ}, \qquad V \in 
\Rep(\cG),
$$
compatible with the tensor structure. Thus, we see that there are 
non-isomorphic objects of $\on{D}(\Gr_G)_\crit\mod$ that go under the 
functor $\Gamma$ to isomorphic objects of $\ghat_\crit\mod_\reg$.

\ssec{}

In \cite{FG2} we showed how to modify the category
$\on{D}(\Gr_G)_\crit\mod$, by simultaneously "adding" to it
$\fZ^\reg_\fg$ as a center, and "dividing" it by the above
$\Rep(\cG)$-action, in order to obtain a category that can be
equivalent to $\ghat_\crit\mod_\reg$.

This procedure amounts to replacing $\on{D}(\Gr_G)_\crit\mod$ by the
appropriate category of {\em Hecke eigen-objects}, denoted
$\on{D}(\Gr_G)_\crit^{\Hecke_\fZ}\mod$.

By definition, an object of $\on{D}(\Gr_G)_\crit^{\Hecke_\fZ}\mod$ is
an object $\F\in \on{D}(\Gr_G)_\crit\mod$, equipped with an action of
the algebra $\fZ^\reg_\fg$ by endomorphisms and a system of
isomorphisms
$$
\al_V: \CF\star \CF_V \overset{\sim}\longrightarrow
\CV_\fZ\underset{\fZ^\reg_\fg}\otimes \CF, \qquad V \in 
\Rep(\cG),
$$
compatible with the tensor structure.

\medskip

We claim that the functor $\Gamma: \on{D}(\Gr_G)_\crit\mod \to
\ghat_\crit\mod_\reg$ naturally gives rise to a functor 
$\Gamma^{\Hecke_\fZ}:\on{D}(\Gr_G)_\crit^{\Hecke_\fZ}\mod\to
\ghat_\crit\mod_\reg$.

\medskip

This is in fact a general property. Suppose for simplicity that we
have an abelian category ${\mc C}$ which is acted upon by the tensor
category $\Rep(H)$, where $H$ is an algebraic group; we denote this
functor by $\CF \mapsto \CF \star V, V \in \Rep(H)$. Let ${\mc
C}^\Hecke$ be the category whose objects are collections
$(\CF,\{\al_V\}_{V \in \Rep(H)})$, where $\CF \in {\mc C}$ and $\{
\al_V \}$ is a compatible system of isomorphisms
$$
\al_V: \CF \star V \overset{\sim}\longrightarrow \underline{V}
\underset{\C}\otimes \CF, \qquad V \in \Rep(H),
$$
where $\underline{V}$ is the vector space underlying $V$. One may
think of ${\mc C}^\Hecke$ as the ``de-equivariantized'' category ${\mc C}$
with respect to the action of $H$. It carries a natural
action of the group $H$: for $h \in H$, we have $h \cdot
(\CF,\{\al_V\}_{V \in \Rep(H)}) = (\CF,\{ (h\otimes \on{id}_\CF) \circ
\al_V\}_{V \in \Rep(H)})$. The category ${\mc C}$ may be reconstructed
as the category of $H$-equivariant objects of ${\mc C}^\Hecke$ with
respect to this action, see \cite{Ga}.

Suppose that we have a functor $\sG: {\mc C} \to {\mc C}'$, such that we
have functorial isomorphisms
\begin{equation}    \label{syst}
\sG(\CF \star V) \simeq \sG(\CF) \underset{\C}\otimes \underline{V},
\qquad V \in \Rep(H),
\end{equation}
compatible with the tensor structure. Then, according to \cite{AG},
there exists a functor ${\sG}^\Hecke: \CC^\Hecke\to \CC'$ such that
$\sG \simeq \sG^\Hecke \circ \on{Ind}$, where the functor $\on{Ind}:
{\mc C} \to {\mc C}^\Hecke$ sends $\CF$ to $\CF \star {\mc O}_{H}$,
where ${\mc O}_{H}$ is the regular representation of $H$. The functor
$\sG^\Hecke$ may be explicitly described as follows: the isomorphisms
$\al_V$ and \eqref{syst} give rise to an action of the algebra ${\mc
O}_{H}$ on $\sG(\CF)$, and ${\sG}^\Hecke(\CF)$ is obtained by taking
the fiber of $\sG(\CF)$ at $1 \in H$.

We take $\CC=\on{D}(\Gr_G)_\crit\mod$, $\CC'=\hg_\crit\mod_\reg$, and
$\sG=\Gamma$. The only difference is that now we are working over the
base $\fZ^\reg_\fg$, which we have to take into account.

\ssec{}

The conjecture suggested in \cite{FG2} states that the resulting functor
\begin{equation} \label{glob Hecke sections}
\Gamma^{\Hecke_\fZ}:\on{D}(\Gr_G)_\crit^{\Hecke_\fZ}\mod\to
\ghat_\crit\mod_\reg.
\end{equation}
is an equivalence. In {\it loc. cit.} we have shown that the functor
$\Gamma^{\Hecke_\fZ}$, when extended to the derived category,
is fully faithful.

\medskip

This conjecture has a number of interesting corollaries pertaining to
the structure of the category of representations at the critical level:

Let us fix a point $\chi\in \Spec(\fZ^\reg_\fg)$, and let us choose a
trivialization of the fiber $\CP_{\cG,\chi}$ of $\CP_{\cG,\fZ}$ at
$\chi$. Let $\ghat_\crit\mod_\chi$ be the subcategory of
$\hg_\crit\mod$, consisting of objects, on which the center acts
according to the character corresponding to $\chi$.

Let $\on{D}(\Gr_G)_\crit^{\Hecke}\mod$ be the category, obtained from
$\on{D}(\Gr_G)_\crit\mod$, by the procedure $\CC\mapsto \CC^\Hecke$
for $H=\cG$, described above. Our conjecture implies that we have an
equivalence
\begin{equation} \label{sect with char}
\on{D}(\Gr_G)_\crit^{\Hecke}\mod\simeq \ghat_\crit\mod_\chi.
\end{equation}
In particular, we obtain that for every two points $\chi,\chi'\in
\Spec(\fZ^\reg_\fg)$ and an isomorphism of $\cG$-torsors
$\CP_{\cG,\chi}\simeq \CP_{\cG,\chi'}$ there exists a canonical
equivalence $\ghat_\crit\mod_\chi\simeq \ghat_\crit\mod_{\chi'}$. This
may be viewed as an analogue of the translation principle that
compares the subcategories $\fg\mod_{\chi}\subset \fg\mod$ for various
central characters $\chi\in \Spec(Z(\fg))$ in the finite-dimensional
case.

By taking $\chi=\chi'$, we obtain that the group $\cG$, or, rather, its
twist with respect to $\CP_{\cG,\chi}$, acts on $\ghat_\crit\mod_\chi$.

\medskip

As we explained in the Introduction to \cite{FG2}, the conjectural
equivalence of \eqref{sect with char} fits into the general picture of
local geometric Langlands correspondence.

Namely, for a point $\chi\in \Spec(\fZ^\reg_\fg)\simeq \Op_\cg(\D)$ as
above, both sides of the equivalence \eqref{sect with char} are
natural candidates for the conjectural Langlands category associated
to the trivial $\cG$-local system on the disc ${\mc D}$. This
category, equipped with an action of the loop group $G\ppart$, should
be thought of as a "categorification" of an irreducible unramified
representation of the group $G$ over a local non-archimedian
field. Proving this conjecture would therefore be the first test of
the local geometric Langlands correspondence proposed in \cite{FG2}.

\ssec{}

Unfortunately, at the moment we are unable to prove the equivalence
\eqref{glob Hecke sections} in general. In this paper we will treat
the following particular case:

Recall that $I^0$ denotes the unipotent radical of the Iwahori
subgroup, and let us consider the corresponding $I^0$-equivariant 
subcategories on both sides of \eqref{glob Hecke sections}.

On the D-module side, we obtain the category
$\on{D}(\Gr_G)_\crit^{\Hecke_\fZ}\mod^{I^0}$, defined in the same way
as $\on{D}(\Gr_G)_\crit^{\Hecke_\fZ}\mod$, but with the requirement
that the underlying D-module $\CF$ be strongly $I^0$-equivariant.

On the representation side, we obtain the category
$\ghat_\crit\mod_\reg^{I^0}$, corresponding to the condition that the
action of $\on{Lie}(I^0)\subset \hg_\crit$ integrates to an algebraic
action of $I^0$.

We shall prove that the functor $\Gamma^{\Hecke_\fZ}$ defines an 
equivalence of categories
\begin{equation}   \label{glob Hecke sections Iw}
\on{D}(\Gr_G)_\crit^{\Hecke_\fZ}\mod^{I^0}\to
\ghat_\crit\mod_\reg^{I^0}.
\end{equation}

This equivalence implies an equivalence 
\begin{equation}    \label{eq for chi}
\on{D}(\Gr_G)_\crit^{\Hecke}\mod^{I^0} \simeq
\ghat_\crit\mod_\chi^{I^0}
\end{equation}
for any fixed character $\chi\in \Spec(\fZ^\reg_\fg)$ and a
trivialization of $\CP_{\cG,\chi}$ as above. In particular, we obtain
the corollaries concerning the translation principle and the action of
$\cG$ on $\ghat_\crit\mod_\chi^{I^0}$.

We remark that from the point of view of the local geometric Langlands
correspondence the categories appearing in the equivalence \eqref{eq
for chi} should be viewed as "categorifications" of the space of
$I$-invariant vectors in an irreducible unramified representation of
the group $G$ over a local non-archimedian field (which is a module
over the corresponding affine Hecke algebra).

\medskip

Let us briefly describe the strategy of the proof. Due to the fact
\cite{FG2} that the functor in one direction in \eqref{glob Hecke
sections Iw} is fully-faithful at the level of the derived categories,
the statement of the theorem is essentially equivalent to the fact
that for every object $\CM\in \ghat_\crit\mod_\reg^{I^0}$ there exists
an object $\CF\in \on{D}(\Gr_G)_\crit^{\Hecke_\fZ}\mod^{I^0}$ and a
non-zero map $\Gamma^{\Hecke_\fZ}(\Gr_G,\CF)\to \CM$. We explain this
in detail in \secref{main proof}.

We exhibit a collection of objects $\BM_{w,\reg}$, numbered by
elements $w\in W$, where $W$ is the Weyl group of $\fg$, which are
quotients of Verma modules over $\ghat_\crit$, such that for every
$\CM\in \ghat_\crit\mod_\reg^{I^0}$ we have
$\on{Hom}(\BM_{w,\reg},\CM)\neq 0$ for at least one $w \in W$.

We then show (see \thmref{get Wakimoto}) that each such $\BM_{w,\reg}$
is isomorphic to $\Gamma^{\Hecke_\fZ}(\Gr_G,\CF_w^\fZ)$ for some
explicit object $\CF_w^\fZ\in
\on{D}(\Gr_G)_\crit^{\Hecke_\fZ}\mod^{I^0}$, thereby proving the
equivalence \eqref{glob Hecke sections Iw}.

\ssec{} 

It is instructive to put our results in the context of other
closely related equivalences of categories.

Using the (tautological) equivalence:
$$\on{D}(\Gr_G)\mod^{I^0}\simeq \on{D}(\wt\Fl_G)\mod^{G[[t]]}$$ (here
and below we omit the subscript $\ka$ when $\ka=0$) and the
equivalence of \thmref{KTthm}, we obtain that for every negative
integral level $\kappa=k\cdot \kappa_{\can}$ there exists an
equivalence between $\on{D}(\Gr_G)\mod^{I^0}$ and the regular block of
the category $\hg_\kappa\mod^{G[[t]]}$, studied in \cite{KL}.  The
latter category is equivalent, according to {\it loc. cit.}, to the
category of modules over the quantum group $U_q^{\on{res}}(\g)$, where
$q = \exp \pi i/(k+h^\vee)$.

Using these equivalences, it was shown in \cite{AG} that the category
$\on{D}(\Gr_G)^{\Hecke}\mod^{I^0}$, defined as above, is equivalent to
the regular block $u_q(\fg)\mod_0$ of the category of modules over the
small quantum group $u_q(\g)$. The tensor product by the line bundle
$\Ll_{\can}^{-h^\vee}$ defines an equivalence
$$\on{D}(\Gr_G)^{\Hecke}\mod^{I^0}\to
\on{D}(\Gr_G)_\crit^{\Hecke}\mod^{I^0}$$ (but this equivalence does
not, of course, respect the functor of global sections). Combining
this with the equivalence of \eqref{eq for chi}, we obtain the
following diagram of equivalent categories:
\begin{equation} \label{quantum group equiv}
\ghat_\crit\mod^{I^0}_\chi \overset{\sim}\leftarrow
\on{D}(\Gr_G)_\crit^{\Hecke}\mod^{I^0} \overset{\sim}\to
u_q(\g)\mod_0.
\end{equation}

\medskip

Recall in addition that in \cite{ABBGM} it was shown that the category
$\on{D}(\Gr_G)^{\Hecke}\mod^{I^0}$ is equivalent to an appropriately
defined category $\on{D}({\mc F}l^{\frac{\infty}{2}})^{I^0}$ of
$I^0$-equivariant D-modules on the semi-infinite flag variety (it is
defined in terms of the Drinfeld compactification
$\ol{\on{Bun}}_N$). Hence, we obtain another diagram of equivalent
categories:
\begin{equation} \label{semiinf equiv}
\ghat_\crit\mod^{I^0}_\chi \overset{\sim}\leftarrow
\on{D}(\Gr_G)_\crit^{\Hecke}\mod^{I^0} \overset{\sim}\to \on{D}({\mc
F}l^{\frac{\infty}{2}})^{I^0}.
\end{equation}

In particular, we obtain a functor
$$\on{D}({\mc F}l^{\frac{\infty}{2}})^{I^0}\to
\ghat_\crit\mod^{I^0}_\chi,$$ which is, moreover, an equivalence. Its
existence had been predicted by B. Feigin and the first named author.

In fact, one would like to be able to define the category $\on{D}({\mc
F}l^{\frac{\infty}{2}})$ without imposing the $I^0$-equivariance
condition, and extend the equivalence of \cite{ABBGM} to this more
general context. Together with the equivalence of \eqref{glob Hecke
sections}, this would imply the existence of the diagram
$$\ghat_\crit\mod_\chi \overset{\sim}\leftarrow
\on{D}(\Gr_G)_\crit^{\Hecke}\mod \overset{\sim}\to \on{D}({\mc
F}l^{\frac{\infty}{2}}),$$ but we are far from being able to achieve
this goal at present.

\medskip

Finally, let us mention one more closely related category, namely, the
derived category $D\bigl(\QCoh((\cG/\cB)^{DG}\mod)\bigr)$ of complexes
of quasi-coherent sheaves over the DG-scheme
$$(\cG/\cB)^{DG}: =
\Spec\Bigl(\Sym_{\CO_{\cG/\cB}}(\Omega^1(\cG/\cB)[1])\Bigr).$$ The
above DG-scheme can be realized as the derived Cartesian product
$\wt\cg\underset{\cg}\times \on{pt}$, where $\on{pt}\to \cg$
corresponds to the point $0\in \cg$, and $\wt\cg = \{ (x,\check\bb)|x
\in \check\bb \subset \cg \}$ is Grothendieck's
alteration.

{}From the results of \cite{ABG} one can obtain an equivalence of the
derived categories
$$
D^b\bigl(\QCoh((\cG/\cB)^{DG}\mod)\bigr) \simeq
D^b\bigl(\on{D}(\Gr_G)^{\Hecke}\mod\bigr)^{I^0}.
$$
Hence we obtain an equivalence:
\begin{equation} \label{coherent equiv}
D^b\bigl(\QCoh((\cG/\cB)^{DG}\mod)\bigr)\simeq
D^b\bigl(\ghat_\crit\mod_\chi\bigr)^{I^0}.
\end{equation}

The existence of such an equivalence follows from the Main Conjecture
6.11 of \cite{FG2}. Note that, unlike the other equivalences mentioned
above, it does not preserve the t-structures, and so is inherently an
equivalence of derived categories.

\ssec{Contents}

Let us briefly describe how this paper is organized:

\medskip

In \secref{Hecke ctry}, after recalling some previous results, we
state the main result of this paper, \thmref{main}.  In \secref{corol}
we review representation-theoretic corollaries of \thmref{main}. In
\secref{main proof} we show how to derive \thmref{main} from
\thmref{get Wakimoto}, and in \secref{sect get Wakimoto} we prove
\thmref{get Wakimoto}.

Finally, in the Appendix, we prove a partial localization result at
the negative level mentioned in \secref{neg level}.

\medskip

The notation in this paper follows that of \cite{FG2}.

\section{The Hecke category}   \label{Hecke ctry}

In this section we recall the main definitions and state our main
result. We will rely on the concepts introduced in our previous paper
\cite{FG2}.

\ssec{Recollections}    \label{recol}

Let $\fg$ be a simple finite-dimensional Lie algebra, and $G$ the
connected algebraic group of adjoint type with the Lie algebra
$\fg$. We shall fix a Borel subgroup $B\subset G$.
Let $\cG$ denote the Langlands dual group of $G$, and by
$\cg$ its Lie algebra. 

\medskip

Let $\Gr_G = G\ppart/G[[t]]$ be the affine Grassmannian associated to
$G$. We denote by $\on{D}(\Gr_G)_\crit\mod$ the category of critically
twisted right D-modules on the affine Grassmannian and by
$\on{D}(\Gr_G)_\crit\mod^{G[[t]]}$ the corresponding
$G[[t]]$-equivariant category. Recall that via the geometric Satake
equivalence (see \cite{MV}) the category
$\on{D}(\Gr_G)_\crit\mod^{G[[t]]}$ has a natural structure of tensor
category under convolution, and as such it is equivalent to
$\Rep(\cG)$.  We shall denote by $V\mapsto \CF_V$ the corresponding
tensor functor $\Rep(\cG) \to \on{D}(\Gr_G)_\crit\mod^{G[[t]]}$.

We have the convolution product functors
$$
\CF \in \on{D}(\Gr_G)_\crit\mod, \CF_V \in
\on{D}(\Gr_G)_\crit\mod^{G[[t]]} \mapsto \CF \star \CF_V \in
\on{D}(\Gr_G)_\crit\mod.
$$
These functors define an action of $\Rep(\cG)$, on the category 
$\on{D}(\Gr_G)_\crit\mod$. Thus, in the terminology of \cite{Ga},
$\on{D}(\Gr_G)_\crit\mod^{G[[t]]}$ has the structure of category
over the stack $\on{pt}/\cG$.

\medskip

Now let $\hg_\crit$ be the {\em affine Kac-Moody algebra} associated
to the critical inner product $-h^\vee \kappa_{\can}$ and
$\hg_\crit\mod$ the category of discrete $\hg_\crit$-modules (see
\cite{FG2}). Its objects are $\hg_\crit$-modules in which every vector
is annihilated by the Lie subalgebra $\fg \otimes t^n\BC[[t]]$ for
sufficiently large $n$. Let $\BV_\crit\in \hg_\crit\mod$ be the vacuum
module $\on{Ind}_{\fg[[t]]\oplus \BC}^{\hg_\crit}(\BC)$. Denote by
$\fZ_\fg$ the topological commutative algebra that is the center of
$\hg_\crit\mod$. Let $\fZ^\reg_\fg$ denote its "regular" quotient,
i.e., the quotient modulo the annihilator of $\BV_\crit$.  We denote
by $\hg_\crit\mod_\reg$ the full subcategory of $\hg_\crit\mod$,
consisting of objects, on which the action of the center $\fZ_\fg$
factors through $\fZ^\reg_\fg$.

\medskip

Recall that via the Feigin-Frenkel isomorphism \cite{FF,F:wak}, the
algebra $\fZ^\reg_\fg$ is identified with the algebra of regular
functions on the scheme $\Op_{\cg}(\D)$ of $\cg$-opers on the formal
disc $\D$. In particular, $\Spec(\fZ^\reg_\fg)$ carries a canonical
$\cG$-torsor, denoted $\CP_{\cG,\fZ}$, whose fiber $\CP_{\cG,\chi}$ at
$\chi \in \Spec(\fZ^\reg_\fg) \simeq \Op_{\cg}(\D)$ is the fiber of
the $\cG$-torsor underlying the oper $\chi$ at the origin of the disc
$\D$. The $\cG$-torsor $\CP_{\cG,\fZ}$ gives rise to a morphism
$\Spec(\fZ^\reg_\fg)\to \on{pt}/\cG$.  We shall denote by $$V\mapsto
\CV_\fZ$$ the resulting tensor functor from $\Rep(\cG)$ to the
category of locally free $\fZ^\reg_\fg$-modules.

\medskip

We define $\on{D}(\Gr_G)_\crit^{\Hecke_\fZ}\mod$ as the fiber product
category
$$\on{D}(\Gr_G)_\crit\mod\underset{\on{pt}/\cG}\times
\Spec(\fZ^\reg_\fg),$$ in the terminology of \cite{Ga}.

\medskip

Explicitly, $\on{D}(\Gr_G)_\crit^{\Hecke_\fZ}\mod$ has as objects the
data of $(\CF,\alpha_V,\,\, \forall\,\, V\in \Rep(\cG))$, where $\CF$
is an object of $\on{D}(\Gr_G)_\crit\mod$, endowed with an action of
the algebra $\fZ^\reg_\fg$ by endomorphisms, and $\alpha_V$ are
isomorphisms of D-modules
$$\CF\star \CF_V\simeq \CV_\fZ\underset{\fZ^\reg_\fg}\otimes \CF,$$
compatible with the action of $\fZ^\reg_\fg$ on both sides, and such
that the following two conditions are satisfied:

\begin{itemize}

\item
For $V$ being the trivial representations $\BC$, the morphism $\alpha_V$
is the identity map.

\item
For $V,W\in \Rep(\cG)$ and $U:=V\otimes W$, the diagram
$$
\CD
(\CF\star \CF_V)\star \CF_W  @>{\sim}>>  \CF\star \CF_U \\
@V{\alpha_V\star \on{id}_{\CF_W}}VV  @V{\alpha_U}VV \\
(\CV_\fZ\underset{\fZ^\reg_\fg}\otimes \CF)\star \CF_W  & &
\CU_\fZ\underset{\fZ^\reg_\fg}\otimes \CF \\
@V{\sim}VV   @V{\sim}VV  \\
\CV_\fZ\underset{\fZ^\reg_\fg}\otimes (\CF\star \CF_W) 
@>{\on{id}_{\CV_\fZ}\otimes \alpha_W}>>  
\CV_\fZ\underset{\fZ^\reg_\fg} \otimes
\CW_\fZ\underset{\fZ^\reg_\fg} \otimes\CF
\endCD
$$ is commutative.

\end{itemize}

Morphisms in this category between $(\CF,\alpha_V)$ and $(\CF',\alpha'_V)$ 
are maps of D-modules $\phi:\CF\to \CF'$ that are compatible with the actions
of $\fZ^\reg_\fg$ on both sides, and such that
$$(\on{id}_{\CV_\fZ}\otimes \phi)\circ \alpha_V=
\alpha'_V\circ (\phi\star \on{id}_{\CF_V}).$$

\ssec{Definition of the functor}    \label{d of f}

Recall that according to \cite{FG1}, the functor of global sections
$$\CF\mapsto \Gamma(\Gr_G,\CF)$$ defines an exact and faithful functor
$\on{D}(\Gr_G)_\crit\mod\to \hg_\crit\mod_\reg$.  Let us recall,
following \cite{FG2}, the construction of the functor
$$\Gamma^{\Hecke_\fZ}:\on{D}(\Gr_G)^{\Hecke_\fZ}_\crit\to
\hg_\crit\mod_\reg.$$

First, let us recall the following result of \cite{BD} (combined with
an observation of \cite{FG2}, Lemma 8.4.3):

\begin{thm}  \label{BDisom} \hfill

\smallskip

\noindent{\em(1)}
For $\CF\in \on{D}(\Gr_G)_\crit\mod$ and $V\in \Rep(\cG)$
we have a functorial isomorphism
$$\beta_V:\Gamma(\Gr_G,\CF\star \CF_V)\simeq 
\Gamma(\Gr_G,\CF)\underset{\fZ^\reg_\fg} \otimes \CV_\fZ.$$

\noindent{\em(2)}
For $\CF,V$ as above and $W\in \Rep(\cG)$, $U:=V\otimes W$ the diagram
$$
\CD
\Gamma(\Gr_G,(\CF\star \CF_V)\star \CF_W) @>{\sim}>> 
\Gamma(\Gr_G,\CF\star (\CF_V\star \CF_W)) \\
@V{\beta_W}VV       @V{\sim}VV   \\
\Gamma(\Gr_G,(\CF\star \CF_V))
\underset{\fZ^\reg_\fg} \otimes \CW_\fZ & & \Gamma(\Gr_G,\CF\star
\CF_U) \\
@V{\beta_V}VV  @V{\beta_U}VV  \\
\Gamma(\Gr_G,\CF)\underset{\fZ^\reg_\fg} \otimes \CV_\fZ
\underset{\fZ^\reg_\fg} \otimes \CW_\fZ @>{\sim}>>
\Gamma(\Gr_G,\CF)\underset{\fZ^\reg_\fg} \otimes \CU_\fZ
\endCD
$$
is commutative.
\end{thm}

\medskip

Consider the scheme
$\Isom_{\fZ}:\Spec(\fZ^\reg_\fg\underset{\on{pt}/\cG}\times
\fZ^\reg_\fg)$.  Let ${{\bf 1}_{\Isom_{\fZ}}}$ denote the unit section
$\Spec(\fZ^\reg_\fg)\to \Isom_{\fZ}$.

Let us denote by
$R_{\fZ}$ the direct image of the structure sheaf under
$\Spec(\fZ^\reg_\fg)\to \on{pt}/\cG$, viewed as an object of
$\Rep(\cG)$. It carries an action of $\fZ^\reg_\fg$ by endomorphisms.
Let $\CR_{\fZ}$ be the associated (infinite-dimensional) vector
bundle over $\Spec(\fZ^\reg_\fg)$; by definition, we have a
canonical isomorphism
$$\CR_\fZ\simeq \Fun(\Isom_\fZ).$$ 
We will think of the projection $p_r:\Isom_\fZ\to \Spec(\fZ^\reg_\fg)$ as 
corresponding to the original $\fZ^\reg_\fg$-action on $R_\fZ$, and hence
on $\CR_\fZ$, by the transport of structure. We will think of the other
projection $p_l:\Isom_\fZ\to \Spec(\fZ^\reg_\fg)$, as corresponding
to the $\fZ^\reg_\fg$-module structure on $\CR_\fZ$ coming from
the fact that this is a vector bundle associated to a $\cG$-representation.

\medskip

We claim that for every object $\CF^H\in
\on{D}(\Gr_G)^{\Hecke_\fZ}_\crit\mod$, the $\hg_\crit$-module
$\Gamma(\Gr_G,\CF^H)$ carries a natural action of the algebra
$\Fun(\Isom_{\fZ})$ by endomorphisms. 

First, note that
$\Gamma(\Gr_G,\CF^H)$ is a $\fZ^\reg_\fg$-bimodule:
we shall refer to the $\fZ^\reg_\fg$-action coming from its
action on any object of $\hg_\crit\mod_\reg$ as "right",
and to the one. coming from the $\fZ^\reg_\fg$-action
on $\CF^H$ as "left".

On the one hand, we have:
$$\Gamma(\Gr_G, \CF^H\star \CF_{R_{\fZ}})\overset{\beta_{R_{\fZ}}}
\simeq \Gamma(\Gr_G, \CF^H)\underset{r,\fZ^\reg_\fg,l}
\otimes\Fun(\Isom_{\fZ}),$$ and on the other hand,
$$\Gamma(\Gr_G, \CF^H\star \CF_{R_{\fZ}})\overset{\alpha_{R_{\fZ}}}
\simeq \Fun(\Isom_{\fZ})\underset{l,\fZ^\reg_\fg,l}\otimes
\Gamma(\Gr_G, \CF^H)\otimes\Fun(\Isom_{\fZ}).$$

By composing we obtain the desired action map
$$\Gamma(\Gr_G, \CF)\underset{r,\fZ^\reg_\fg,l}\otimes \Fun(\Isom_{\fZ})
\overset{\alpha_{R_\fZ}\circ \beta^{-1}_{R_\fZ}}\longrightarrow
\Fun(\Isom_{\fZ})\underset{l,\fZ^\reg_\fg,l}\otimes \Gamma(\Gr_G, \CF^H)
\overset{{\bf 1}^*_{\Isom_{\fZ}}}\longrightarrow \Gamma(\Gr_G, \CF^H).$$
The fact that it is associative follows from the second condition on
$\alpha_V$ and \thmref{BDisom}(2).

\medskip

We define the functor $\Gamma^{\Hecke_\fZ}$ by
$$\CF^H\mapsto \Gamma(\Gr_G,\CF^H) \underset{\Fun(\Isom_{\fZ}),{\bf
1}^*_{\Isom_{\fZ}}}\otimes \fZ^\reg_\fg.$$ Since the functor $\Gamma$
is exact, the functor $\Gamma^{\Hecke_\fZ}$ is evidently right-exact,
and we will denote by $\on{L}\Gamma^{\Hecke_\fZ}$ its left derived
functor $D^-(\on{D}(\Gr_G)^{\Hecke_\fZ}_\crit\mod)\to
D^-(\hg_\crit\mod_\reg)$

\medskip

The following was established in \cite{FG2}, Theorem 8.7.1:

\begin{thm}  \label{GH fully faithful}
The functor $\on{L}\Gamma^{\Hecke_\fZ}$, restricted to
$D^b(\on{D}(\Gr_G)^{\Hecke_\fZ}_\crit\mod)$, is fully faithful.
\end{thm}

In \cite{FG2} we formulated the following

\begin{conj}    \label{general conj}
The functor $\Gamma^{\Hecke_\fZ}$ is exact and defines an equivalence
of categories $\on{D}(\Gr_G)^{\Hecke_\fZ}_\crit\mod$ and
$\hg_\crit\mod_\reg$.
\end{conj}

\ssec{The statement of the main result}    \label{main result}

Recall that both categories $\on{D}(\Gr_G)^{\Hecke_\fZ}_\crit\mod$ and
$\hg_\crit\mod_\reg$ carry a natural action of the group $G\ppart$
(see \cite{FG2}, Sect. 22, where this is discussed in detail). Let
$I\subset G[[t]]$ be the Iwahori subgroup, the preimage of the Borel
subgroup $B \subset G$ in $G[[t]]$ under the evaluation map $G[[t]]
\to G$. Let $I^0$ be the unipotent radical of $I$. Let us denote by
$\on{D}(\Gr_G)^{\Hecke_\fZ}_\crit\mod^{I^0}$ and
$\hg_\crit\mod_\reg^{I^0}$ the corresponding categories if
$I$-equivariant objects. Since $I^0$ is connected, these are full
subcategories in $\on{D}(\Gr_G)^{\Hecke_\fZ}_\crit\mod$ and
$\hg_\crit\mod_\reg^{I^0}$, respectively.


The functor $\Gamma^{\Hecke_\fZ}$ induces a functor
$\on{D}(\Gr_G)^{\Hecke_\fZ}_\crit\mod^{I^0}\to
\hg_\crit\mod_\reg^{I^0}$. 
The goal of the present paper is to prove the following:

\begin{thm}  \label{main} \hfill

\smallskip

\noindent{\em (1)} 
For any $\CF^H\in \on{D}(\Gr_G)^{\Hecke_\fZ}_\crit\mod^{I^0}$ we have
$L^i \Gamma^{\Hecke_\fZ}(\Gr_G,\CF^H) = 0$ for all $i>0$.

\smallskip

\noindent{\em (2)} The functor $$\Gamma^{\Hecke_\fZ}:
\on{D}(\Gr_G)^{\Hecke_\fZ}_\crit\mod^{I^0}\to
\hg_\crit\mod_\reg^{I^0}$$ is an equivalence of categories.
\end{thm}

This is a special case of \conjref{general conj}.

\section{Corollaries of the main theorem}   \label{corol}

We shall now discuss some applications of \thmref{main}. Note that
both sides of the equivalence stated in \thmref{main} are categories
over the algebra $\fZ^\reg_\fg$.

\ssec{Specialization to a fixed central character} \label{disc of cor}

Let us fix a point $\chi\in\Spec(\fZ^\reg_\fg)$, i.e., a character of
$\fZ^\reg_\fg$, and consider the subcategories on both sides of the
equivalence of \thmref{main}(2), corresponding to objects on which the
center acts according to this character.  Let us denote the resulting
subcategory of $\hg_\crit\mod_\reg^{I^0}$ by
$\hg_\crit\mod_\chi^{I^0}$. The resulting subcategory of
$\on{D}(\Gr_G)^{\Hecke_\fZ}_\crit\mod^{I^0}$ can be described as
follows.

\medskip

Let us denote by $\on{D}(\Gr_G)^{\Hecke}_\crit\mod$ the category,
whose objects are the data of $(\CF,\alpha_V)$, where $\CF\in
\on{D}(\Gr_G)_\crit\mod$ and $\alpha_V$ are isomorphisms of D-modules
defined for every $V\in \Rep(\cG)$
$$\CF\star\CF_V\simeq \uV\underset{\BC}\otimes \CF,$$ where $\uV$
denotes the vector space underlying the representation $V$. These
isomorphisms must be compatible with tensor products of objects of
$\Rep(\cG)$ in the same sense as in the definition of
$\on{D}(\Gr_G)^{\Hecke_\fZ}_\crit\mod$.

Note that $\on{D}(\Gr_G)^{\Hecke}_\crit\mod$ carries a natural weak
action of the algebraic group $\cG$: \footnote{We refer the reader to
\cite{FG2}, Sect. 20.1, where this notion is introduced.} Given an
$S$-point $\bg$ of $\cG$ and an $S$-family of objects $(\CF,\alpha_V)$
of $\on{D}(\Gr_G)^{\Hecke}_\crit\mod$ we obtain a new $S$-family by
keeping $\CF$ the same, but replacing $\alpha_V$ by $\bg\cdot
\alpha_V$, where $\bg$ acts naturally on $\uV\otimes \CO_S$.

In addition, $\on{D}(\Gr_G)^{\Hecke}_\crit\mod$ carries a commuting
Harish-Chandra action of the group $G\ppart$; in particular, the
subcategory $\on{D}(\Gr_G)^{\Hecke}_\crit\mod^{I^0}$ makes sense.

\medskip

Let $\CP_{\cG,\chi}$ be the fiber of the $\cG$-torsor $\CP_{\cG,\fZ}$
at $\chi$. Tautologically we have:

\begin{lem} \hfill

\smallskip

\noindent{\em (1)} For every trivialization
$\gamma:\CP_{\cG,\chi}\simeq \CP^0_\cG$ there exists a canonical
equivalence respecting the action of $G\ppart$
$$\bigl(\on{D}(\Gr_G)^{\Hecke_\fZ}_\crit\mod\bigr)_\chi\simeq
\on{D}(\Gr_G)^{\Hecke}_\crit\mod,$$ where the LHS denotes the fiber of
$\on{D}(\Gr_G)^{\Hecke_\fZ}_\crit\mod$ at $\chi$.

\smallskip

\noindent{\em (2)} If $\gamma'=\bg\cdot \gamma$ for $\bg\in \cG$, the
above equivalence is modified by the self-functor of
$\on{D}(\Gr_G)^{\Hecke}_\crit\mod$, given by the action of $\bg$.

\end{lem}

Hence, from \thmref{main} we obtain:

\begin{cor}    \label{main cor}
For every trivialization $\gamma:\CP_{\cG,\chi}\simeq \CP^0_\cG$ there
exists a canonical equivalence
$$\hg_\crit\mod_\chi^{I^0}\simeq
\on{D}(\Gr_G)^{\Hecke}_\crit\mod^{I^0}.$$
\end{cor}


\medskip

{}From \corref{main cor} we obtain:

\begin{cor} \hfill

\smallskip

\noindent{\em (1)} For any two points $\chi_1,\chi_2\in
\Spec(\fZ^\reg_\fg)$ and an isomorphism of $\cG$-torsors
$\CP_{\cG,\chi_1}\simeq \CP_{\cG,\chi_2}$ there exists a canonical
equivalence
$$\hg_\crit\mod_{\chi_1}^{I^0}\simeq \hg_\crit\mod_{\chi_2}^{I^0}.$$

\smallskip

\noindent{\em (2)} For every $\chi\in \Spec(\fZ^\reg_\fg)$, the group
of automorphisms of the $\cG$-torsor $\CP_{\cG,\chi}$ acts on the
category $\hg_\crit\mod_\chi^{I^0}$.

\end{cor}

More generally, let $S$ be an affine scheme, and let $\chi_{1,S}$ and 
$\chi_{2,S}$ be two $S$-points of $\Spec(\fZ^\reg_\fg)$. Let
$\hg_\crit\mod_{S,1}^{I^0}$ and $\hg_\crit\mod_{S,2}^{I^0}$ be the
corresponding base-changed categories. 

By definition, the
objects of $\hg_\crit\mod_{i,S}$ are the objects of $\hg_\crit\mod_\reg$,
endowed with an action of $\CO_S$ compatible with the initial action
of $\fZ^\reg_\fg$ on $\CM$ via the homomorphism
$\fZ^\reg_\fg\to \CO_S$, corresponding to $\chi_{i,S}$.
Morphisms in this category are $\hg_\crit$-morphisms compatible with 
the action of $\CO_S$. 

\medskip

We obtain:

\begin{cor} \label{main cor, families}
For every lift of the map 
$$(\chi_{1,S}\times \chi_{2,S}):S\to \Spec(\fZ^\reg_\fg)\times
\Spec(\fZ^\reg_\fg)$$ to a map $S\to \Isom_\fZ$, there exists a
canonical equivalence
$$\hg_\crit\mod_{S,1}^{I^0}\simeq \hg_\crit\mod_{S,2}^{I^0}.$$
\end{cor}

\ssec{Description of irreducibles}   \label{sect descr of irr}

\corref{main cor} allows to describe explicitly the set of irreducible
objects in $\hg_\crit\mod_\reg^{I^0}$. For that we will need to recall
some more notation related to the categories
$\on{D}(\Gr_G)^{\Hecke}_\crit\mod$ and
$\on{D}(\Gr_G)^{\Hecke_\fZ}_\crit\mod$.

\medskip

Consider the forgetful functor $\on{D}(\Gr_G)^{\Hecke}_\crit\mod\to
\on{D}(\Gr_G)_\crit\mod$. It admits a left adjoint, denoted 
$\on{Ind}^\Hecke$, which can be described as follows. 

Let $R$ be the object of $\Rep(\cG)$ equal to $\CO_\cG$ under the
left regular action; let $\CF_R$ denote the corresponding
object of $\on{D}(\Gr_G)_\crit\mod^{G[[t]]}$. Then for $\CF\in
\on{D}(\Gr_G)_\crit\mod$, the convolution $\CF\star \CF_R$ is naturally
an object of $\on{D}(\Gr_G)^{\Hecke}_\crit\mod$, and it is easy to see
that $\on{Ind}^\Hecke(\CF):=\CF\star \CF_R$ is the desired left adjoint.

\medskip

Similarly, the forgetful functor
$\on{D}(\Gr_G)^{\Hecke_\fZ}_\crit\mod\to \on{D}(\Gr_G)_\crit\mod$
admits a left adjoint functor $\on{Ind}^{\Hecke_\fZ}$ given by
$\CF\mapsto \CF\star \CF_{R_\fZ}$. The next assertion follows from the
definitions:

\begin{lem}  \label{H induced}  \hfill

\smallskip

\noindent{\em (1)} For $\CF\in \on{D}(\Gr_G)_\crit\mod$ there exist
canonical isomorphisms:
$$\Gamma(\Gr_G,\on{Ind}^{\Hecke_\fZ}(\CF))\simeq
\Gamma(\Gr,\CF)\underset{\fZ^\reg_\fg}\otimes \Fun(\Isom_{\fZ}),$$
where $\Fun(\Isom_{\fZ})$ is a module over $\fZ^\reg_\fg$ via one
of the projections $\Isom_{\fZ}\to \Spec(\fZ^\reg_\fg)$.

\smallskip

\noindent{\em (2)}
For $\CF$ as above,
$$\Gamma^{\Hecke_\fZ}\bigl(\Gr_G,
\on{Ind}^{\Hecke_\fZ}(\CF)\bigr)\simeq \Gamma(\Gr,\CF).$$
\end{lem}

\medskip

Let us now recall the description of irreducible objects of
$\on{D}(\Gr_G)^{\Hecke}_\crit\mod^{I^0}$, established in \cite{ABBGM},
Corollary 1.3.10.

Recall that $I$-orbits on $\Gr_G$ are parameterized by the set
$W_{\aff}/W$, where $W_{\aff}$ denotes the extended affine Weyl
group. For an element $\wt{w}\in W_{\aff}$ let us denote by
$\IC_{\wt{w},\Gr_G}$ the corresponding irreducible object of
$\on{D}(\Gr_G)_\crit\mod^I$.

For an element $w\in W$, let $\cla_w\in W_{\aff}$ denote the unique
dominant coweight satisfying:
$$
\begin{cases}
&\langle \alpha_\imath,\cla\rangle=0 \text{ if } w(\alpha_\imath) \text{ is
positive, and} \\ &\langle \alpha_\imath,\cla\rangle=1 \text{ if }
w(\alpha_\imath) \text{ is negative,}
\end{cases}
$$
for $\imath$ running over the set of vertices of the Dynkin diagram. 

It was shown in {\it loc. cit.} that the objects
$\Ind^\Hecke(\IC_{w\cdot \lambda_w})$ for $w\in W$ are the
irreducibles of $\on{D}(\Gr_G)^{\Hecke}_\crit\mod^{I^0}$.

\medskip

Combining this with \lemref{H induced} and \corref{main cor}, we
obtain:

\begin{thm} \label{decsr of irr}
Isomorphism classes of irreducible objects of
$\hg_\crit\mod^{I^0}_\reg$ are parameterized by pairs $(\chi\in
\Spec(\fZ^\reg_\fg),w\in W)$. For each such pair the corresponding
irreducible object is given by
$$\Gamma(\Gr_G,\IC_{w\cdot \lambda_w})\underset{\fZ^\reg_\fg}\otimes
\BC_\chi.$$
\end{thm}

\ssec{The algebroid action}   \label{algebroid}

Let $\isom_{\fZ}$ be the Lie algebroid of the groupoid
$\Isom_{\fZ}$. According to \cite{BD} (see also \cite{FG2}, Sect. 7.4
for a review), we have a canonical action of $\isom_{\fZ}$ on
$\wt{U}^\reg_\crit(\hg)$ by outer derivations, where
$\wt{U}^\reg_\crit(\hg)$ is the topological associative algebra
corresponding to the category $\hg_\crit\mod_\reg$ and its
tautological forgetful functor to vector spaces.  

In more detail, there exists a topological associative algebra,
denoted $U^{\ren,\reg}(\hg_\crit)$, and called 
the renormalized universal enveloping algebra at the critical
level. It is endowed with a natural filtration, with the $0$-th term
$U^{\ren,\reg}(\hg_\crit)_0$ being $U^\reg(\hg_\crit)$, and 
$$U^{\ren,\reg}(\hg_\crit)_1/U^{\ren,\reg}(\hg_\crit)_0\simeq 
U^\reg(\hg_\crit)\underset{\fZ^\reg_\fg}\hattimes \isom_\fZ.$$

The action of $\isom_{\fZ}$ on $\wt{U}^\reg_\crit(\hg)$ is given
by the adjoint action of  
$\isom_\fZ$, regarded as a subset of
$\subset U^{\ren,\reg}(\hg_\crit)_1/U^{\ren,\reg}(\hg_\crit)_0$.

\medskip

Let $S$ be an affine scheme, and let $\chi_S$ be an $S$-point of
$\Spec(\fZ^\reg_\fg)$.  Let $\xi_S$ be a section of
$\isom_{\fZ}|_S$. Set $S':=S\times \Spec(\BC[\epsilon]/\epsilon^2)$;
then the image of $\xi_S$ in $T(\Spec(\fZ^\reg_\fg))|_S$ gives rise to
an $S'$-point, denoted, $\chi'_S$, of $\Spec(\fZ^\reg_\fg)$.

Let $\hg_\crit\mod_S$ (resp., $\hg_\crit\mod_{S'}$) be the
corresponding base-changed category, where the latter identifies with
the category of discrete modules over
$\wt{U}^\reg_\crit(\hg)\underset{\fZ^\reg_\fg}\otimes \CO_S$ (resp.,
$\wt{U}^\reg_\crit(\hg)\underset{\fZ^\reg_\fg}\otimes \CO_{S'}$).
Then the above action of $\isom_{\fZ}$ on $\hg_\crit\mod_\reg$ gives
rise to the following construction:

\medskip

To every $\CM\in \hg_\crit\mod_S$ we can functorially attach an
extension
\begin{equation}  \label{action of algebroid} 
0\to \CM\to \CM'\to \CM\to 0,\,\, \qquad \CM'\in \hg_\crit\mod_{S'}.
\end{equation}

The module $\CM'$ is defined as follows. The above action of
$\isom_{\fZ}$ by outer derivations of $\wt{U}^\reg_\crit(\hg)$ allows
to lift $\xi_S$ to an isomorphism
$$A(\xi_S):\wt{U}^\reg_\crit(\hg)\underset{\fZ^\reg_\fg,\chi_{S'}}
\otimes \CO_{S'}\to
\wt{U}^\reg_\crit(\hg)\underset{\fZ^\reg_\fg,\chi_S} \otimes
\CO_S[\epsilon]/\epsilon^2.$$ We set $\CM'$ to be the
$\wt{U}^\reg_\crit(\hg)\underset{\fZ^\reg_\fg,\chi_{S'}} \otimes
\CO_{S'}$-module, corresponding via $A(\xi_S)$ to
$\CM[\epsilon]/\epsilon^2$.

The isomorphism $A(\xi_S)$ is defined up to conjugation by an element
of the form $1+\epsilon \cdot u$, $u\in
\wt{U}^\reg_\crit(\hg)\underset{\fZ^\reg_\fg,\chi_S} \otimes
\CO_S$. Since this automorphism can be canonically lifted onto
$\CM[\epsilon]/\epsilon^2$, we obtain that $\CM'$ is well-defined.

By construction, the functor $\CM\mapsto \CM'$ respects the
Harish-Chandra $G\ppart$-actions on the categories $\hg_\crit\mod_{S}$
and $\hg_\crit\mod_{S'}$, respectively.

\medskip

Let us note now that a data $(\chi_S:S\to \Spec(\fZ^\reg_\fg),
\xi_S\in \isom_{\fZ}|_S)$ as above can be regarded as a map $S'\to
\Isom_{\fZ}$, where first and second projections
$$S'\to \Isom_{\fZ}\rightrightarrows \Spec(\fZ^\reg_\fg)$$ 
are equal to
$$S'\to S\overset{\chi_S}\to \Spec(\fZ^\reg_\fg) \text{ and }
S'\overset{\chi'_S}\to \Spec(\fZ^\reg_\fg),$$
respectively.

Hence, \corref{main cor, families} gives rise to an equivalence
$$\hg_\crit\mod^{I^0}_S\otimes \BC[\epsilon]/\epsilon^2\simeq
\hg_\crit\mod^{I^0}_{S'},$$ and, in particular, to a functor

\begin{equation} \label{second action of algebroid}
\hg_\crit\mod^{I^0}_S\to \hg_\crit\mod^{I^0}_{S'}.
\end{equation}

\begin{prop} \label{two actions of algebroid}
The functor 
$$\CM\mapsto \CM':\hg_\crit\mod_S\to \hg_\crit\mod_{S'}$$ of
\eqref{action of algebroid}, restricted to $\hg_\crit\mod^{I^0}_S$, is
canonically isomorphic to the above functor \eqref{second action of
algebroid}. 
\end{prop}

\begin{proof}

The assertion follows from the fact that for $\CF\in
\on{D}(\Gr_G)_\crit\mod$, the $\hg_\crit$-action on
$\Gamma(\Gr_G,\CF)$ lifts canonically to an action of
$U^{\ren,\reg}(\hg_\crit)$ (see \cite{FG2}, Sect 7.4), so that for
$(S,\chi_S,\xi_S)$ as above we have a canonical trivialization
$$\gamma_\CF:\Gamma(\Gr_G,\CF)'\simeq
\Gamma(\Gr_G,\CF)[\epsilon]/\epsilon^2,$$ in the notation of
\eqref{action of algebroid}. Moreover, this functorial isomorphism is
compatible with that of \thmref{BDisom} in the sense that for every
$V\in \Rep(\cG)$, the diagram
$$
\CD \Gamma(\Gr_G,\CF\star \CF_V)' @>{\gamma_{\CF\star \CF_V}}>>
\Gamma(\Gr_G,\CF\star \CF_V)[\epsilon]/\epsilon^2 \\ @V{\beta_V}VV
@V{\beta_V\otimes \on{id}}VV \\
\Bigl(\Gamma(\Gr_G,\CF)\underset{\fZ^\reg_\fg}\otimes \CV\Bigr)'
@>{\gamma_\CF\otimes \xi_S}>>
\Bigl(\Gamma(\Gr_G,\CF)\underset{\fZ^\reg_\fg}\otimes
\CV\Bigr)[\epsilon]/\epsilon^2, \endCD
$$ commutes, where the bottom arrow comprises the isomorphism
$\gamma_\CF$ and the canonical action of $\xi_S$ on $\CV_\fZ$. The
latter compatibility follows assertion (b) in Theorem 8.4.2 of
\cite{FG2}.

\end{proof}

\ssec{Relation to semi-infinite cohomology}

Let us consider the functor of semi-infinite cohomology on the category
$\hg_\crit\mod_{\reg}^{I^0}$:
$$\CM \mapsto H^\semiinfb(\fn\ppart,\fn[[t]],\CM \otimes \Psi_0)$$
(see \cite{FG2}, Sect. 18 for details concerning this functor).

For an $S$-point $\chi_S$ of $\Spec(\fZ^\reg_\fg)$ and $\CM\in
\hg_\crit\mod_S$, each $H^\semiinfi(\fn\ppart,\fn[[t]],\CM\otimes
\Psi_0)$ is naturally an $\CO_S$-module.

\medskip

Let now $(\chi_{1,S},\chi_{2,S})$ be a pair of $S$-points of
$\Spec(\fZ^\reg_\fg)$, equipped with a lift $S\to \Isom_\fZ$,
and let $\CM_1\in \hg_\crit\mod_{S,1}^{I^0}$
and $\CM_2\in \hg_\crit\mod_{S,2}^{I^0}$ be two objects corresponding
to each other under the equivalence of \corref{main cor, families}.

\begin{prop}  \label{behaviour of semiinf}
Under the above circumstances the $\CO_S$-modules
$$H^\semiinfi(\fn\ppart,\fn[[t]],\CM_1\otimes \Psi_0)\text{ and }
H^\semiinfi(\fn\ppart,\fn[[t]],\CM_2\otimes \Psi_0)$$
are canonically isomorphic.
\end{prop}

\begin{proof}

The assertion of the proposition can be tautologically translated as
follows:

The functor
$$\on{D}(\Gr_G)_\crit\mod \overset{\Gamma} \to \hg_\crit\mod_\reg
\overset{H^\semiinfi(\fn\ppart,\fn[[t]],?\otimes
\Psi_0)}\longrightarrow \fZ^\reg_\fg\mod$$ factors through a functor
$$H^\semiinfi_{\cG}:\on{D}(\Gr_G)_\crit\mod \to \Rep(\cG),$$ followed
by the pull-back functor, corresponding to the morphism
$\Spec(\fZ^\reg_\fg)\to \on{pt}/\cG$.  Moreover, for $V\in \Rep(\cG)$
we have a functorial isomorphism
\begin{equation} \label{semiinf ident}
H^\semiinfi_{\cG}(\CF\star \CF_V)\simeq H^\semiinfi_{\cG}(\CF)\otimes V,
\end{equation}
compatible with the isomorphism of \thmref{BDisom}(1).

\medskip

The sought-after functor $H^\semiinfi_{\cG}$ has been
essentially constructed in \cite{FG2}, Sect. 18.3. Namely,
$$\Hom_{\cG}\bigl(V^\cla,H^\semiinfi_{\cG}(\CF)\bigr):=
H^i(N\ppart,\CF|_{N\ppart\cdot t^\cla}\otimes \Psi_0),$$ in the
notation of {\it loc. cit.} The isomorphisms \eqref{semiinf ident}
follow from the definitions.

\end{proof}

Finally, we would like to compare the isomorphisms of
\propref{behaviour of semiinf} and \propref{two actions of
algebroid}. Let $\CM$ be an object of $\hg_\crit\mod_\reg^{I^0}$; let
$\chi_S$ be an $S$-point of $\Spec(\fZ^\reg_\fg)$ and $\xi_S$ a
section of $\isom_\fZ|_S$.

On the one hand, in Proposition 18.3.2 of \cite{FG2} we have shown
that there exists a canonical isomorphism:
$${\bf a}_\CM:H^\semiinfi(\fn\ppart,\fn[[t]],\CM'\otimes \Psi_0)\simeq
H^\semiinfi(\fn\ppart,\fn[[t]],\CM\otimes \Psi_0)[\epsilon]/\epsilon^2,$$
valid for any $\CM\in \hg_\crit\mod_\reg$. 

On the other hand, combining \propref{two actions of algebroid} and 
\propref{behaviour of semiinf} we obtain another isomorphism
$${\bf b}_\CM:H^\semiinfi(\fn\ppart,\fn[[t]],\CM'\otimes \Psi_0)\simeq
H^\semiinfi(\fn\ppart,\fn[[t]],\CM\otimes \Psi_0)[\epsilon]/\epsilon^2.$$

Unraveling the two constructions, we obtain the following:

\begin{lem}
The isomorphisms ${\bf a}_\CM$ and ${\bf b}_\CM$ coincide.
\end{lem}

\section{Proof of the main theorem}   \label{main proof}

In \secref{main result} we have constructed a functor
$$
\Gamma^{\Hecke_\fZ}:
\on{D}(\Gr_G)^{\Hecke_\fZ}_\crit\mod^{I^0}\to
\hg_\crit\mod_\reg^{I^0}.
$$
Now we wish to show that this functor is an equivalence of
categories. This will prove \thmref{main}.

We start by constructing in \secref{Fw} certain objects $\CF^\fZ_w, w
\in W$, of the category $\on{D}(\Gr_G)^{\Hecke_\fZ}_\crit\mod^{I^0}$
such that $\Gamma^{\Hecke_\fZ}(\CF^\fZ_w) \simeq {\mathbb
M}_{w,\on{reg}}$, the ``standard modules'' of the category
$\hg_\crit\mod_\reg^{I^0}$. The main result of \secref{Fw},
\thmref{get Wakimoto}, will be proved in \secref{sect get
Wakimoto}. Next, in \secref{part one} we prove part (1) of
\thmref{main} that the functor $\Gamma^{\Hecke_\fZ}$ is exact. We then
outline in \secref{general} a general framework for proving that it is
an equivalence. Using this framework, we prove \thmref{main} modulo
\thmref{get Wakimoto}.

In \secref{remark on general} we explain what needs to be done in
order to prove our stronger \conjref{general conj}. Finally, in
Sects.~\ref{another}--\ref{proof of flatness} we give an alternative
proof of part (1) of \thmref{main}.

\ssec{Standard modules}    \label{Fw}

For an element $w\in W$, let $\BM_w$ be the Verma module over $\ghat$,
$$
\BM_w = \on{Ind}^{\hg_\crit}_{\fg[[t]]}(M_{w(\rho)-\rho}),
$$
where for a weight $\lambda$ we denote by $M_\lambda$ the Verma module
over $\fg$ with highest weight $\lambda$.

Let $\BM_{w,\reg}$ be the maximal quotient module that belongs to
$\hg_\crit\mod_\reg$, i.e.,
$\BM_{w,\reg}=\BM_{w}\underset{\fZ_\fg}\otimes \fZ_\fg^\reg$. In fact,
it was shown in \cite{FG2}, Corollary 13.3.2, that as modules over
$\fZ_\fg$, all $\BM_w$ are supported over a quotient algebra
$\fZ^\nilp_\fg$, and are flat as $\fZ^\nilp_\fg$-modules. The
subscheme $\Spec(\fZ^\reg_\fg)\subset \Spec(\fZ_\fg)$ is contained in
$\Spec(\fZ^\nilp_\fg)$, so the definition of $\BM_{w,\reg}$ does not
neglect any lower cohomology.

The main ingredient in the remaining steps of our proof of
\thmref{main} is the following:

\begin{thm}  \label{get Wakimoto}
For each $w\in W$ there exists an 
object $\CF^\fZ_w\in \on{D}(\Gr_G)^{\Hecke_\fZ}_\crit\mod^{I^0}$, such that
$\Gamma^{\Hecke_\fZ}(\Gr_G, \CF_w)$ is isomorphic to $\BM_{w,\reg}$.
\end{thm}

The proof of this theorem will consist of an explicit construction of
the objects $\CF^\fZ_w$, which will be carried out in \secref{sect get
Wakimoto}.

The proof of \thmref{main} will only use a part of the assertion of
\thmref{get Wakimoto}: namely, that there exist objects $\CF^\fZ_w\in
\on{D}(\Gr_G)^{\Hecke_\fZ}_\crit\mod^{I^0}$, endowed with a surjection
\begin{equation} \label{surj on Wak} 
\Gamma^{\Hecke_\fZ}(\Gr_G, \CF^\fZ_w)\twoheadrightarrow \BM_{w,\reg}.
\end{equation}
What we will actually use is the following corollary of this
statement:

\begin{cor} \label{F non-zero}
For every $\CM\in \hg_\crit\mod_\reg^{I^0}$ there exists
an object $\CF^H\in \on{D}(\Gr_G)^{\Hecke_\fZ}_\crit\mod^{I^0}$
and a non-zero map $\Gamma^{\Hecke_\fZ}(\Gr_G, \CF^H)\to \CM$.
\end{cor}

\begin{proof}

By \cite{FG2}, Lemma 7.8.1, for every object $\CM\in
\hg_\crit\mod_\reg^{I^0}$ there exists $w\in W$ and a non-zero map
$\BM_{w,\reg}\to \CM$.

\end{proof}

\ssec{Exactness}    \label{part one}

Let us recall from \secref{sect descr of irr} the left adjoint functor
$\on{Ind}^{\Hecke_\fZ}$ to the obvious forgetful functor
$\on{D}(\Gr_G)^{\Hecke_\fZ}_\crit\mod \to \on{D}(\Gr_G)_\crit\mod$.

It is clear that every object of
$\on{D}(\Gr_G)^{\Hecke_\fZ}_\crit\mod$ can be covered by one of the
form $\on{Ind}^{\Hecke_\fZ}(\CF)$. From \lemref{H induced}(1) we
obtain that we can use bounded from above complexes, whose terms
consist of objects of the form $\on{Ind}^{\Hecke_\fZ}(\CF)$, in order
to compute $\on{L}\Gamma^{\Hecke_\fZ}$. Thus, we obtain:

\begin{lem}  \label{L as tor}
For $\CF^H\in \on{D}(\Gr_G)^{\Hecke_\fZ}_\crit\mod$,
$$\on{L}^i\Gamma^{\Hecke_\fZ}(\Gr_G,\CF^H)\simeq
\on{Tor}_i^{\Fun(\Isom_{\fZ})}\bigl(\Gamma(\Gr_G,\CF^H),
\fZ^\reg_\fg\bigr).$$
\end{lem}

\medskip

We shall call an object of $\on{D}(\Gr_G)^{\Hecke_\fZ}_\crit\mod$ {\it
finitely generated} if it can be obtained as a quotient of an object
of the form $\on{Ind}^{\Hecke_\fZ}(\CF)$, where $\CF$ is a finitely
generated object of $\on{D}(\Gr_G)_\crit\mod$.

It is easy to see that an object $\CF^H\in
\on{D}(\Gr_G)^{\Hecke_\fZ}_\crit\mod$ is finitely generated if and
only if the functor
$\Hom_{\on{D}(\Gr_G)^{\Hecke_\fZ}_\crit\mod}(\CF^H,\cdot)$ commutes
with direct sums.

We shall call an object of $\on{D}(\Gr_G)^{\Hecke_\fZ}_\crit\mod$ {\it
finitely presented}, if it is isomorphic to
$\on{coker}\bigl(\on{Ind}^{\Hecke_\fZ}(\CF_1)\to
\on{Ind}^{\Hecke_\fZ}(\CF_2)\bigr)$, where $\CF_1,\CF_2$ are both
finitely generated objects of $\on{D}(\Gr_G)_\crit\mod$. The following
lemma is straightforward.

\begin{lem}   \label{all as ind fp} \hfill

\smallskip

\noindent{\em (1)} An object $\CF^H\in
\on{D}(\Gr_G)^{\Hecke_\fZ}_\crit\mod$ is finitely presented if and
only if the functor
$\Hom_{\on{D}(\Gr_G)^{\Hecke_\fZ}_\crit\mod}(\CF^H,\cdot)$ commutes
with filtering direct limits.

\smallskip

\noindent{\em (2)} Every object of
$\on{D}(\Gr_G)^{\Hecke_\fZ}_\crit\mod$ is isomorphic to a filtering
direct limit of finitely presented ones.

\end{lem}

The proof of the following proposition will be given in \secref{sect
finite amplitude}.

\begin{prop}   \label{finite amplitude}
For every finitely presented object of
$\on{D}(\Gr_G)^{\Hecke_\fZ}_\crit\mod$, the corresponding object
$\on{L}\Gamma^{\Hecke_\fZ}(\Gr_G,\CF^H)\in D^-(\hg_\crit\mod_\reg)$
belongs to $D^b(\hg_\crit\mod_\reg)$.
\end{prop}

The crucial step in the proof of part (1) of \thmref{main} is the
following:

\begin{prop}  \label{bounded exact}
If $\CF^H\in \on{D}(\Gr_G)^{\Hecke_\fZ}_\crit\mod^{I^0}$ is such that
$\on{L}\Gamma^{\Hecke_\fZ}(\Gr_G,\CF^H)$ belongs to 
$D^b(\hg_\crit\mod_\reg)^{I^0}$, then 
$$\on{L}^i\Gamma^{\Hecke_\fZ}(\Gr_G,\CF^H)\bigr)=0, \qquad i>0.$$
\end{prop}

\begin{proof}
Let $\CM$ be the lowest cohomology of
$\on{L}\Gamma^{\Hecke_\fZ}(\Gr_G,\CF^H)$, which lives, say, in degree
$-k$. By \corref{F non-zero} there exists another object $\CF_1^H\in
\on{D}(\Gr_G)^{\Hecke_\fZ}_\crit\mod^{I^0}$ and a non-zero map
$\Gamma^{\Hecke_\fZ}(\Gr_G,\CF^H_1)\to \CM$. Hence, we obtain a
non-zero map in $D^-(\hg_\crit\mod_\reg)$
$$\on{L}\Gamma^{\Hecke_\fZ}(\Gr_G,\CF_1^H)[k]\to
\on{L}\Gamma^{\Hecke_\fZ}(\Gr_G,\CF^H).$$

But by \thmref{GH fully faithful}, such map comes from a map
$\CF_1^H[k]\to \CF^H$, which is impossible if $k>0$.
\end{proof}

\medskip

\noindent{\em Proof of part (1) of \thmref{main}}.  Combining
\propref{finite amplitude} and \propref{bounded exact}, we obtain that
$\on{L}^i\Gamma^{\Hecke_\fZ}(\Gr_G,\CF^H)=0$ for any $i>0$ and any
$\CF^H\in \on{D}(\Gr_G)^{\Hecke_\fZ}_\crit\mod^{I^0}$, which is
finitely presented.

However, by \lemref{L as tor}, the functors
$$\CF^H\mapsto \on{L}^i\Gamma^{\Hecke_\fZ}(\Gr_G,\CF^H)$$
commute with direct limits, and our assertion
follows from \lemref{all as ind fp}(2).\qed

\ssec{Proof of the equivalence}    \label{general}

Consider the following general categorical framework. Let
$\sG:\CC_1\to \CC_2$ be an exact functor between abelian
categories. Assume that for $X,Y\in \CC_1$ the maps
$$\Hom_{\CC_1}(X,Y)\to \Hom_{\CC_2}(\sG(X),\sG(Y)) \text{ and }
\Ext^1_{\CC_1}(X,Y)\to \Ext^1_{\CC_2}(\sG(X),\sG(Y))$$
are isomorphisms. 

\begin{lem}    \label{adjoint}
If $\sG$ admits a right adjoint functor $\sF$ which is conservative,
then $\sG$ is an equivalence. \footnote{Recall that a functor $\sF$ is
called {\em conservative} if for any $X\neq 0$ we have $\sF(X)\neq
0$.}
\end{lem}

\begin{proof}

The fully faithfulness assumption on $\sG$ implies that the adjunction
map induces an isomorphism between the composition $\sF\circ \sG$ and
the identity functor on $\CC_1$.  We have to show that the second
adjunction map is also an isomorphism.

For $X'\in \CC_2$ let $Y'$ and $Z'$ be the kernel and cokernel,
respectively, of the adjunction map
$$\sG\circ \sF(X')\to X'.$$

Being a right adjoint functor, $\sF$ is left-exact, hence we obtain an
exact sequence
$$0\to \sF(Y')\to \sF\circ \sG\circ \sF(X')\to \sF(X').$$ But since
$\sF(X')\to \sF\circ \sG(\sF(X'))$ is an isomorphism, we obtain that
$\sF(Y')=0$. Since $\sF$ is conservative, this implies that $Y'=0$.

\medskip

Suppose that $Z'\neq 0$. Since $\sF(Z')\neq 0$, there exists an object
$Z\in \CC_1$ with a non-zero map $\sG(Z)\to Z'$. Consider the induced
extension
$$0\to \sG\circ \sF(X')\to W'\to \sG(Z)\to 0.$$ Since $\sG$ induces a
bijection on $\Ext^1$, this extension can be obtained from an
extension
$$0\to \sF(X')\to W\to Z\to 0$$ in $\CC_1$. In other words, we obtain
a map $\sG(W)\to X'$, which does not factor through $\sG\circ
\sF(X')\subset X'$, which contradicts the $(\sG,\sF)$ adjunction.

\end{proof}

Thus, in order to prove of part (2) of \thmref{main} it remains to show
that the functor
$\Gamma^{\Hecke_\fZ}:\on{D}(\Gr_G)^{\Hecke_\fZ}_\crit\mod^{I^0}\to
\hg_\crit\mod_\reg^{I^0}$ admits a right adjoint. (The fact that it is
conservative will then follow immediately from \corref{F non-zero}.)

Recall from \cite{FG2}, Sect. 20.7, that the tautological functor
$\on{D}(\Gr_G)^{\Hecke_\fZ}_\crit\mod^{I^0}\hookrightarrow
\on{D}(\Gr_G)^{\Hecke_\fZ}_\crit\mod$ admits a right adjoint, given by
$\on{Av}_{I^0}$.  Hence, it suffices to prove the following:

\begin{prop}   \label{exists right adjoint}
The functor $\Gamma^{\Hecke_\fZ}:\on{D}(\Gr_G)^{\Hecke_\fZ}_\crit\mod\to
\hg_\crit\mod_\reg$ admits a right adjoint.
\end{prop}

\begin{proof}

First, we will show the following:

\begin{lem}  \label{kappa adjoint}
The functor $\Gamma:\on{D}(\Gr_G)_\crit\mod\to \hg_\crit\mod_\reg$
admits a right adjoint.
\end{lem}

\begin{proof}

We will prove that for any level $k$ the functor
$\Gamma:\on{D}(\Gr_G)_k\mod\to \hg_\crit\mod_k$ admits a
right adjoint (see the Introduction for the definition of these
categories). I.e., we have to prove the representability of the
functor
\begin{equation} \label{which functor}
\CF\mapsto \Hom_{\hg_k\mod}\bigl(\Gamma(\Gr_G,\CF),\CM\bigr)
\end{equation}
for every given $\CM\in \hg_k\mod$.

\medskip

Consider the following general set-up. Let $\CC$ be an abelian
category, and let $\CC^0$ be a full (but not necessarily abelian)
subcategory, such that the following holds:

\begin{itemize}

\item $\CC^0$ is equivalent to a small category.

\item The cokernel of any surjection $X''\twoheadrightarrow X'$ with
$X',X''\in \CC^0$, also belongs to $\CC^0$.

\item $\CC$ is closed under filtering direct limits.

\item For $X\in \CC^0$, the functor $\Hom_{\CC}(X,\cdot)$ commutes 
with filtering direct limits.

\item Every object of $\CC$ is isomorphic to a filtering direct limit
of objects of $\CC^0$.

\end{itemize}

Then we claim that any contravariant left-exact functor $\sF\to
\Vect$, which maps direct sums to direct products (and, hence, direct
limits to inverse limits, by the previous assumption), is
representable.

Indeed, given such $\sF$, consider the category of pairs $(X,f)$,
where $X\in \CC^0$ and $f\in \sF(X)$. Morphisms between $(X,f)$ and
$(X',f')$ are maps $\phi:X\to X'$, such that $\phi^*(f')=f$. By the
first assumption on $\CC^0$, this category is small. By the second
assumption on $\CC^0$ and the left-exactness of $\sF$, this category
is filtering. It is easy to see that the object
$$\underset{\underset{(X,f)}\longrightarrow}{\lim}\, X.$$
represents the functor $\sF$.

\medskip

We apply this lemma to $\CC=\on{D}(\Gr_G)_k\mod$ with $\CC^0$
being the subcategory of finitely-generated D-modules. We set $\sF$ to
be the functor \eqref{which functor}, and the representability
assertion follows.

Note that we could have applied the above general principle to
$\CC=\on{D}(\Gr_G)^{\Hecke_\fZ}_\crit\mod$ and $\CC^0$ being the
subcategory of finitely presented objects, and obtain the assertion
of \propref{exists right adjoint} right away.

\end{proof}

Thus, for $\CM$, let $\CF$ be the object of $\on{D}(\Gr_G)_\crit\mod$
that represents the functor
$$\CF_1\mapsto
\Hom_{\hg_\crit\mod_\reg}\bigl(\Gamma(\Gr_G,\CF_1),\CM\bigr)$$ for a
given $\CM\in \hg_\crit\mod_\reg$. We claim that $\CF$ is naturally an
object of $\on{D}(\Gr_G)^{\Hecke_\fZ}_\crit\mod$ and that it
represents the functor
\begin{equation} \label{which functor two}
\CF^H_1\mapsto 
\Hom_{\hg_\crit\mod_\reg}\bigl(\Gamma^{\Hecke_\fZ}(\Gr_G,\CF^H_1),\CM\bigr).
\end{equation}

First, since the algebra $\fZ^\reg_\fg$ acts on $\CM$ by
endomorphisms, the object $\CF$ carries an action of $\fZ^\reg_\fg$ by
functoriality. Let us now construct the morphisms
$\alpha_V$. Evidently, it is sufficient to do so for $V$
finite-dimensional. Let $V^*$ denote its dual.

For a test object $\CF_1\in \on{D}(\Gr_G)_\crit\mod$ 
we have:
\begin{align*}
&\Hom_{\on{D}(\Gr_G)_\crit\mod}(\CF_1,\CF\star \CF_V)\simeq
\Hom_{\on{D}(\Gr_G)_\crit\mod}(\CF_1\star \CF_{V^*},\CF)\simeq \\
&\simeq \Hom_{\hg_\crit\mod_\reg}\bigl(\Gamma(\Gr_G,\CF_1\star
\CF_{V^*}),\CM\bigr)\simeq \\ &\simeq
\Hom_{\hg_\crit\mod_\reg}\bigl(\Gamma(\Gr_G,\CF_1)
\underset{\fZ^\reg_\fg}\otimes \CV^*_{\fZ^\reg_\fg},\CM\bigr)\simeq
\Hom_{\hg_\crit\mod_\reg}\bigl(\Gamma(\Gr_G,\CF_1),\CV_\fZ
\underset{\fZ^\reg_\fg}\otimes\CM\bigr),
\end{align*}
where the last isomorphism takes place since $\CV_\fZ$ is locally
free.  For the same reason,
$$\Hom_{\on{D}(\Gr_G)_\crit\mod}(\CF_1,\CV_\fZ
\underset{\fZ^\reg_\fg}\otimes\CF)\simeq
\Hom_{\hg_\crit\mod_\reg}\bigl(\Gamma(\Gr_G,\CF_1),\CV_\fZ
\underset{\fZ^\reg_\fg}\otimes\CM\bigr),$$
which implies that there exists a canonical isomorphism $\alpha_V$
$$\CF\star \CF_V\simeq \CV_\fZ
\underset{\fZ^\reg_\fg}\otimes\CF,$$
as required. That these isomorphisms are compatible with tensor
products of objects of $\Rep(\cG)$ follows from \thmref{BDisom}(2).

Finally, the fact that $(\CF,\alpha_V)$, thus defined, represents the
functor \eqref{which functor two}, follows from the construction. This
completes the proof of \propref{exists right adjoint}.

\end{proof}

Thus, we obtain that the functor $\Gamma^{\Hecke_\fZ}$ admits a right
adjoint functor. Moreover, this right adjoint functor is conservative
by \corref{F non-zero}. Therefore part (2) of \thmref{main} now
follows from part (1), proved in \secref{part one}, and
\lemref{adjoint}, modulo \propref{finite amplitude} and \thmref{get
Wakimoto}. It remains to prove those two statements. \propref{finite
amplitude} will be proved in the next subsection and \thmref{get
Wakimoto} will be proved in \secref{sect get Wakimoto}.

\ssec{Proof of \propref{finite amplitude}}   \label{sect finite amplitude}

Recall the category $\on{D}(\Gr_G)^{\Hecke}_\crit\mod$, introduced in
\secref{sect descr of irr}.  Recall also that the $\cG$-torsor
$\CP_{\cG,\fZ}$ on $\Spec(\fZ^\reg_\fg)$ is non-canonically trivial,
and let us fix such a trivialization. This choice identifies the
category $\on{D}(\Gr_G)^{\Hecke_\fZ}_\crit\mod$ with
$\on{D}(\Gr_G)^{\Hecke}_\crit\mod\otimes \fZ^\reg_\fg$, i.e., with the
category of objects of $\on{D}(\Gr_G)^{\Hecke}_\crit\mod$ endowed with
an action of $\fZ^\reg_\fg$ by endomorphisms.

Under this equivalence, the functor $\CF\mapsto
\on{Ind}^{\Hecke_\fZ}(\CF)$ goes over to
$$\CF\mapsto \on{Ind}^{\Hecke}(\CF)\otimes \fZ^\reg_\fg.$$ Note also
that the trivialization of $\CP_{\cG,\fZ}$ identifies $\Isom_{\fZ}$
with $\Spec(\fZ^\reg_\fg)\times \cG\times \Spec(\fZ^\reg_\fg)$, so
that the map ${\bf 1}_{\Isom_{\fZ}}$ corresponds to
$\Delta_{\Spec(\fZ^\reg_\fg)}\times {\bf 1}_{\cG}$.  For $\CF$ as
above, we have an identification
$$\Gamma\bigl(\Gr_G,\on{Ind}^{\Hecke_\fZ}(\CF)\bigr)\simeq
\Gamma(\Gr_G,\CF)\otimes \CO_\cG \otimes \fZ^\reg.$$

\medskip

Let $\CF^H$ be a finitely presented object of
$\on{D}(\Gr_G)^{\Hecke_\fZ}_\crit\mod$ equal to the cokernel of a map
$$\phi:\on{Ind}^{\Hecke}(\CF_1)\otimes \fZ^\reg_\fg\to
\on{Ind}^{\Hecke}(\CF_2)\otimes \fZ^\reg_\fg.$$ Recall that
$\fZ^\reg_\fg$ is isomorphic to a polynomial algebra
$\BC[x_1,...,x_n,...]$.  Since $\CF_1$ was assumed finitely generated,
a map as above has the form $\phi_m\otimes
\on{id}_{\BC[x_{m+1},x_{m+2},...]}$, where $\phi_m$ is a map
$$\on{Ind}^{\Hecke}(\CF_1)\otimes \BC[x_1,...,x_m]\to 
\on{Ind}^{\Hecke}(\CF_2)\otimes \BC[x_1,...,x_m]$$
defined for some $m$.

Hence, as a module over $\Fun(\Isom_{\fZ})\simeq \fZ^\reg_\fg\otimes
\CO_\cG\otimes \fZ^\reg_\fg$,
\begin{equation} \label{shape of module}
\Gamma(\Gr_G,\CF^H)\simeq \CL\otimes \BC[x_{m+1},x_{m+2},...],
\end{equation}
where $\CL$ is some module over $\fZ^\reg_\fg\otimes \CO_\cG\otimes
\BC[x_1,...,x_m]$.

\medskip

We can compute
$$\Gamma(\Gr_G,\CF^H)\underset{\Fun(\Isom_{\fZ})}{\overset{L}\otimes}
\fZ^\reg_\fg$$ in two steps, by first restricting to the preimage of
the diagonal under
$$\Spec(\fZ^\reg_\fg)\times \cG\times
\Spec(\fZ^\reg_\fg)\twoheadrightarrow
\Spec(\BC[x_{m+1},x_{m+2},...])\times
\Spec(\BC[x_{m+1},x_{m+2},...]),$$ and then by further restriction to
$\Spec(\BC[x_1,...,x_m])\times \Spec(\BC[x_{m+1},x_{m+2},...])$
sitting inside
$$\Spec(\BC[x_1,...,x_m])\times \cG\times \Spec(\BC[x_1,...,x_m])\times
\Spec(\BC[x_{m+1},x_{m+2},...]).$$

When we apply the first step to the module appearing in \eqref{shape
of module}, it is acyclic off cohomological degree $0$. The second
step has a cohomological amplitude bounded by $m+\dim(\cG)$.

Hence,
$$\on{Tor}_i^{\Fun(\Isom_{\fZ})}\bigl(\Gamma(\Gr_G,\CF^H),
\fZ^\reg_\fg\bigr)=0$$ for $i>m+\dim(\cG)$, which is what we had to
show.

This completes the proof of \propref{finite amplitude}. Therefore the
proof of \thmref{main} is now complete modulo \thmref{get Wakimoto}.

\ssec{A remark on the general case}    \label{remark on general}

Let us note that the proof of \thmref{main} presented above would
enable us to prove the general \conjref{general conj} if we could
show that the functor
$$\Loc:\hg_\crit\mod_{\reg}\to \on{D}(\Gr_G)_\crit\mod,$$
right adjoint to the functor
$\Gamma:\on{D}(\Gr_G)_\crit\mod\to \hg_\crit\mod_\reg$
is conservative. In other words, in order to prove \conjref{general conj}
we need to know that for every $\CM\in \hg_\crit\mod_\reg$ there
exists a critically twisted D-module $\CF$ on $\Gr_G$ with
a non-zero map $\Gamma(\Gr_G,\CF)\to \CM$. This,
in turn, can be reformulated as follows:

Let $\on{Diff}(\Gr_G)_\crit$ be the *-sheaf of critically twisted
differential operators on $\Gr_G$. This is a pro-object of
$\on{D}(\Gr_G)_\crit\mod$, defined by the property that
$$\on{Hom}(\on{Diff}(\Gr_G)_\crit,\CF)\simeq \Gamma(\Gr_G,\CF)$$
functorially in $\CF\in \on{D}(\Gr_G)_\crit\mod$.

Explicitly, let us write $\Gr_G$ as
$\underset{\CY}{\underset{\longrightarrow}{"\lim"}}\, \CY$, where
$\CY\subset \Gr_G$ are closed sub-schemes. For each such $\CY$, let
$\on{Dist}(\CY)_\crit\in \on{D}(\Gr_G)_\crit\mod$ be the twisted
D-module of distributions on $\CY$, i.e., the object
$\on{Ind}^{\on{D}(\Gr_G)_\crit\mod}_{\QCoh(\Gr_G)}(\CO_{\CY})$, which
means by definition that
$$\on{Hom}_{\on{D}(\Gr_G)_\crit\mod}\Bigl(
\on{Ind}^{\on{D}(\Gr_G)_\crit\mod}_{\QCoh(\Gr_G)}(\CO_{\CY}),\CF\Bigr)=
\on{Hom}_{\QCoh(\Gr_G)}\Bigl(\CO_Y,\CF\Bigr).$$

Then
$$\on{Diff}(\Gr_G)_\crit: =
\underset{\CY}{\underset{\longleftarrow}{"\lim"}}\,
\on{Dist}(\CY)_\crit\in \on{Pro}(\on{D}(\Gr_G)_\crit\mod).$$

Let $\Gamma(\Gr_G,\on{Diff}(\Gr_G)_\crit)$ be the corresponding object
of $\on{Pro}(\hg_\crit\mod_\reg)$.

We obtain:

\begin{cor}
The following assertions are equivalent:

\smallskip

\noindent{\em (1)}
\conjref{general conj} holds.

\smallskip

\noindent{\em (2)}
The object $\Gamma(\Gr_G,\on{Diff}(\Gr_G)_\crit)$ is a pro-projective
generator of $\hg_\crit\mod_\reg$.

\smallskip

\noindent{\em (3)}
The functor on $\hg_\crit\mod_\reg$
$$\CM\mapsto \on{Hom}_{\hg_\crit\mod_\reg}\Bigl(
\Gr_G,\on{Diff}(\Gr_G)_\crit,\CM\Bigr)$$ is conservative.

\end{cor}

\ssec{Another proof of exactness}    \label{another}

In this subsection we give shall present an alternative proof of part
(1) of \thmref{main}.

According to \lemref{L as tor}, proving the exactness property stated
in part (1) of \thmref{main} is equivalent to proving that
\begin{equation} \label{van tor}
\on{Tor}_i^{\Fun(\Isom_{\fZ})}\bigl(\Gamma(\Gr_G,\CF^H),
\fZ^\reg_\fg\bigr)=0
\end{equation}
for all $i>0$ and $\CF^H\in
\on{D}(\Gr_G)^{\Hecke_\fZ}_\crit\mod^{I^0}$. We will derive this from
the following weaker statement:

\begin{prop}  \label{flatness}
For every $\CF\in \on{D}(\Gr_G)_\crit\mod^{I^0}$, the space of
sections $\Gamma(\Gr_G,\CF)$ is flat as a $\fZ^\reg_\fg$-module.
\end{prop}

Note that our general conjecture \eqref{general conj} predicts that
both \eqref{van tor} and the assertion of \propref{flatness} should
hold without the $I^0$-equivariance assumption. However, at the moment
we can neither prove the corresponding generalization of
\propref{flatness} nor derive \eqref{van tor} from it.

\medskip

Let us first show how \propref{flatness} implies \eqref{van tor} on
the $I^0$-equivariant category.

\begin{prop}  \label{analysis of irr}
Every finitely generated object of
$\on{D}(\Gr_G)^{\Hecke_\fZ}_\crit\mod^{I^0}$ admits a finite
filtration, whose subquotients are of the form
\begin{equation} \label{form of irr}
\on{Ind}^{\Hecke_\fZ}(\CF)\underset{\fZ^\reg_\fg}\otimes \CL,
\end{equation}
where $\CL$ is a $\fZ^\reg_\fg$-module.
\end{prop}

Let us deduce \eqref{van tor} from this proposition. 

\begin{proof}

It is enough to show that \eqref{van tor} holds for finitely presented
objects of the category $\on{D}(\Gr_G)^{\Hecke_\fZ}_\crit\mod^{I^0}$. By
\propref{analysis of irr}, we conclude that it is enough to consider
objects of $\on{D}(\Gr_G)^{\Hecke_\fZ}_\crit\mod^{I^0}$ of the form
given by \eqref{form of irr}.

We have:
$$\Gamma\bigl(\Gr_G,\on{Ind}^{\Hecke_\fZ}(\CF)\underset{\fZ^\reg_\fg}
\otimes \CL\bigr) \underset{\Fun(\Isom_{\fZ})}{\overset{L}\otimes}
\fZ^\reg_\fg\simeq
\Gamma\bigl(\Gr_G,\CF)\underset{\fZ^\reg_\fg}{\overset{L}\otimes}
\CL,$$ and the assertion follows from \propref{flatness}.

\end{proof}

Let us now prove \propref{analysis of irr}.

\begin{proof}

Choosing a trivialization of $\CP_{\cG,\fZ}$ as in the previous
subsection, we can identify
$\on{D}(\Gr_G)^{\Hecke_\fZ}_\crit\mod^{I^0}$ with
$\on{D}(\Gr_G)^{\Hecke}_\crit\mod^{I^0}\otimes \fZ^\reg_\fg$.

\medskip

Similarly to the case of $\on{D}(\Gr_G)^{\Hecke_\fZ}_\crit\mod$, we
shall call an object of $\on{D}(\Gr_G)^{\Hecke}_\crit\mod$ finitely
generated if it is isomorphic to a quotient of some
$\on{Ind}^\Hecke(\CF)$ for a finitely generated $\CF\in
\on{D}(\Gr_G)_\crit\mod$.

\medskip

Let us recall from \cite{ABBGM}, Corollary 1.3.10(1), that every
finitely generated object in $\on{D}(\Gr_G)^{\Hecke}_\crit\mod^{I^0}$
has a finite length. Therefore, every finitely generated object of
$\on{D}(\Gr_G)^{\Hecke}_\crit\mod^{I^0}\otimes \fZ^\reg_\fg$ admits a
finite filtration, whose subquotients are quotients of modules of the
form $\CF^H\otimes \fZ^\reg_\fg$ with $\CF^H\in
\on{D}(\Gr_G)^{\Hecke}_\crit\mod^{I^0}$ being irreducible.  However,
every such quotient has the form $\CF^H\otimes \CL$ for some
$\fZ^\reg_\fg$-module $\CL$.

Moreover, as was mentioned in \secref{sect descr of irr}, by
\cite{ABBGM}, Corollary 1.3.10(2), every irreducible in
$\on{D}(\Gr_G)^{\Hecke}_\crit\mod^{I^0}$ is of the form
$\on{Ind}^\Hecke(\CF)$ for some $\CF\in
\on{D}(\Gr_G)_\crit\mod^{I^0}$. This implies the assertion of the
proposition.

\end{proof}

\ssec{Proof of \propref{flatness}}    \label{proof of flatness}

We can assume that our object $\CF\in \on{D}(\Gr_G)_\crit\mod^{I^0}$
is finitely generated, which automatically implies that it has a
finite length. This reduces us to the case when $\CF$ is irreducible.

It is easy to see that any irreducible object of
$\on{D}(\Gr_G)_\crit\mod^{I^0}$ is equivariant also with respect to
$\BG_m$, which acts on $G\ppart$, and hence on $\Gr_G$, by rescalings
$t \mapsto at$. Moreover, the grading arising on its space of sections
is bounded from above. (Our conventions are such that $\BV_\crit$ is
{\it negatively} graded.)

\medskip

Recall now that the action of $\wt{U}^\reg_\crit(\hg)$ on a module of
the form $\Gamma(\Gr_G,\CF)$ for an object $\CF\in \on{D}(\Gr_G)_\crit\mod$
canonically extends to an action of the renormalized algebra
$U^{\ren,\reg}(\hg_\crit)$.  Recall also that
$U^{\ren,\reg}(\hg_\crit)$ contains a $\fZ^\reg_\fg$ sub-bimodule and
a Lie subalgebra $\wt{U}^\reg_\crit(\hg)^\sharp$, which is an
extension
$$0\to \wt{U}^\reg_\crit(\hg)\to \wt{U}^\reg_\crit(\hg)^\sharp\to
\isom_{\fZ}\to 0.$$ (The resulting action of $\isom_{\fZ}$ by outer
derivations on $\wt{U}^\reg_\crit(\hg)$ is the one discussed in
\secref{algebroid}.)

We will prove the following general assertion, which implies
\propref{flatness}:

\begin{lem}
Let $\CM$ be an object of $\hg_\crit\mod_\reg$, such that the action
of $\wt{U}^\reg_\crit(\hg)$ on it extends to an action of
$U^{\ren,\reg}(\hg_\crit)$.  Assume also that $\CM$ is endowed with a
grading, compatible with the one on $U^{\ren,\reg}(\hg_\crit)$, given
by rescalings $t \mapsto at$. Finally, assume that the grading on
$\CM$ is bounded from above. Then $\CM$ is flat as a
$\fZ^\reg_\fg$-module.
\end{lem}

The proof is a variation of the argument used in \cite{BD},
Sect. 6.2.2:

\begin{proof}

We can identify $\fZ^\reg_\fg$ with a polynomial algebra
$\BC[x_1,...,x_n,...]$.  Moreover, we can do so in a
grading-preserving fashion, in which case each generator $x_i$ will be
homogeneous of a {\it negative} degree.

It is enough to show that $\CM$ is flat over each subalgebra 
$\BC[x_1,...,x_m]\subset \fZ^\reg_\fg$. We will prove the following
assertion:

\medskip

\noindent {\it For every vector $\bv\in
\BA^m:=\Spec(\BC[x_1,...,x_m])$, the $\BC[x_1,...,x_m]$-module $\CM$
is (non-canonically) isomorphic to its translate by means of $\bv$.}

\medskip

Clearly, a module over $\BC[x_1,...,x_m]$ having this property is
flat. To prove the above claim we proceed as follows. Choose a section
$\xi$ of $\isom_{\fZ}$, which projects onto $\bv$ under
$\isom_{\fZ}\to T(\Spec(\fZ^\reg_\fg))$, where we think of $\bv$ as a
constant vector field on $\fZ^\reg_\fg\simeq
\Spec(\BC[x_1,...,x_n,...])$. Let us further lift $\xi$ to an element
$\xi'$ of $\wt{U}^\reg_\crit(\hg)^\sharp$.

Since the grading on the $x_i$'s is positive, we can choose $\xi'$
to belong to the (completion of the) sum of strictly positive graded 
components of $\wt{U}^\reg_\crit(\hg)^\sharp$. 

Then the assumption that the grading on $\CM$ is bounded from above,
implies that $\on{exp}(\xi')$ is a well-defined automorphism of $\CM$
as a vector space. This automorphism covers the automorphism
$\on{exp}(\bv)$ of $\BC[x_1,...,x_m]$, and the latter is the same as
the translation by $\bv$.

\end{proof}

\section{Proof of \thmref{get Wakimoto}}  \label{sect get Wakimoto}

In this section we construct the objects $\CF_w^\fZ$ of the category
$\on{D}(\Gr_G)^{\Hecke_\fZ}_\crit\mod^{I^0}$ whose existence is stated
in \thmref{get Wakimoto}.

\ssec{}    \label{wak without center}

We first describe the analogues of these objects in the category
$\on{D}(\Gr_G)^{\Hecke}_\crit\mod^{I^0}$. These objects, which we will
denote by $\CF_w$, were studied in \cite{ABBGM} under the name "baby
co-Verma modules".

\medskip

First, we consider the case $w=w_0$. Recall that the Langlands dual
group comes equipped with a standard Borel subgroup $\cB\subset \cG$;
we shall denote by $\cH$ the Cartan quotient of $\cB$.

Let $\cB^-\subset \cG$ be a Borel subgroup in the generic relative
position with respect to $\cB$. The latter means that $\cB\cap \cB^-$
is {\it a} Cartan subgroup; we shall identify it with $\cH$ by means
of the projection
$$\cB\cap \cB^-\hookrightarrow \cB\twoheadrightarrow \cH.$$

For $\cla\in \cLambda^+$ let $\ell^\cla$ be the line of coinvariants
$(V^\cla)_{\cN^-}$, where $V^\cla$ denotes the standard irreducible 
$\cG$-representation of highest weight $\cla$ with respect to $\cB$. 

The assignment $\cla\mapsto \ell^\cla$ is an $\cH$-torsor, and we obtain
a collection of maps
\begin{equation} \label{simple Plucker}
V^\cla\overset{\kappa^\cla}\twoheadrightarrow \ell^\cla,
\end{equation}
satisfying the Pl\"ucker relations, i.e., for any two dominant
coweights $\cla$ and $\cmu$, the diagram
\begin{equation} \label{Plucker rel}
\CD
V^\cla\otimes V^\cmu @>{\kappa^\cla\otimes \kappa^\cmu}>>
\ell^\cla\otimes \ell^\cmu \\ @VVV @V{\sim}VV \\ V^{\cla+\cmu}
@>{\kappa^{\cla+\cmu}}>> \ell^{\cla+\cmu} \endCD
\end{equation}
commutes.

\medskip

Let $\Fl_G = G\ppart/I$ be the affine flag variety. We have the
category $\on{D}(\Fl_G)_\crit\mod$ of right critically twisted
D-modules on $\Fl_G$ and the corresponding Iwahori equivariant
category $\on{D}(\Fl_G)_\crit\mod^I$. Given $\CF\in
\on{D}(\Gr_G)_\kappa\mod^I$ and $\CM\in \on{D}(\Fl_G)_\crit\mod^I$, we
can form their convolution, denoted by $\CM \underset{I}\star \CF$,
which is an object of $D^b(\on{D}(\Fl_G)_\crit\mod)^I$ (see \cite{FG2}
for details).

\medskip

For a dominant map $\cla$ let $j_{\cla,*}$ denote the $*$-extension of
the critically twisted D-module corresponding to the constant sheaf on
the Iwahori orbit of the point $t^\cla\in \Fl_G$.  Let
$j_{\cla,\Gr_G,*}\in \on{D}(\Gr_G)_\crit\mod^I$ be
$j_{\cla,*}\underset{I}\star \delta_{1,\Gr_G}$; in other words it is
the $*$-extension of the constant D-module on the Iwahori orbit of the
point $t^\cla\in \Gr_G$. Note that for $\cmu\in \cLambda^+$ we have a
canonical map
$$j_{\cla,\Gr_G,*}\star \CF_{V^\cmu}\to  j_{\cla+\cmu,\Gr_G,*},$$
obtained by identifying $\CF_{V^\cmu}$ with $\IC_{\Grb^\cmu}$.

Consider the object of $\on{D}(\Gr_G)^\Hecke_\crit\mod$ equal to the
direct sum
$$\wt{\CF}_{w_0}:=\underset{\cla\in \cLambda^+}\oplus\, 
\on{Ind}^{\Hecke}\bigl(j_{\cla,\Gr_G,*}\bigr) \otimes \ell^{-\cla}.$$

\medskip

For a dominant coweight $\cmu$ we have an evident map
\begin{equation} \label{shift map}
j_{\cmu,*}\underset{I}\star \wt{\CF}_{w_0}\to \ell^\cmu \otimes
\wt{\CF}_{w_0}.
\end{equation} 

We obtain two maps $\wt{\CF}_{w_0}\star \CF_{V^\cmu}\rightrightarrows
\wt{\CF}_{w_0}\otimes \ell^\cmu$ that correspond to the two circuits
of the following non-commutative diagram:
$$
\CD
\wt{\CF}_{w_0}\star \CF_{V^\cmu} @>{\alpha_V}>> \uV^\cmu\otimes
\wt{\CF}_{w_0} \\ @VVV @V{\kappa^\cmu}VV \\
j_{\cmu,*}\underset{I}\star \wt{\CF}_{w_0} @>>> \ell^\cmu \otimes
\wt{\CF}_{w_0}, \endCD
$$
where the left vertical arrow comes from the following map, defined
for each $\cla$:
$$j_{\cla,\Gr_G,*}\star \CF_{R}\star \CF_{V^\cmu}\simeq
j_{\cla,\Gr_G,*}\star \CF_{V^\cmu}\star \CF_R\to
j_{\cla+\cmu,\Gr_G,*}\star \CF_R.$$ Here we are using the object
$\CF_R$ of $\on{D}(\Gr_G)_\crit\mod^{G[[t]]}$ introduced in
\secref{sect descr of irr}, so that $\on{Ind}^\Hecke(\CF)\simeq
\CF\star \CF_R$.

We set $\CF_{w_0}$ to be the maximal quotient of $\wt{\CF}_{w_0}$,
which co-equalizes the resulting two maps
$$\ell^{-\cmu}\otimes \wt{\CF}_{w_0}\star \CF_{V^\cmu}\rightrightarrows
\wt{\CF}_{w_0}$$ for every $\cmu\in \cLambda^+$.
Note that the map \eqref{shift map} gives rise to a map
\begin{equation} \label{shift maps}
j_{\cmu,*}\underset{I}\star \CF_{w_0}\to \ell^\cmu \otimes \CF_{w_0}.
\end{equation}
By construction, $\CF_{w_0}$ has the following universal property:

\medskip

Let $\CF^H$ be an object of $\on{D}(\Gr_G)_\crit\mod^I$, endowed
with a system of morphisms
\begin{equation} \label{shift maps 2}
j_{\cmu,*}\underset{I}\star\CF^H\to \ell^\cmu \otimes \CF^H,
\end{equation}
compatible with the isomorphisms
\begin{equation} \label{j's multiply}
j_{\cmu,*}\underset{I}\star j_{\cmu',*}\simeq j_{\cmu+\cmu',*}
\end{equation}
and $\ell^\cmu\otimes \ell^{\cmu'}\simeq \ell^{\cmu+\cmu'}$.

Let $\phi: \CF_R\to \CF^H$ be a map, such that for every
$\cmu\in \cLambda$ the following diagram is commutative:
$$
\CD
\CF_R\star \CF_{V^\cmu} @>{\alpha_V}>>  \uV^\cmu\otimes R 
@>{\on{id}_{\uV^\cla}\otimes \phi}>> \uV^\cmu\otimes \CF^H 
@>{\kappa^\cmu}>>  \ell^\cmu\otimes \CF^H \\
@V{\sim}VV   & & & & @AAA  \\
\CF_{V^\cmu}\star \CF_R @>>> j_{\cmu,\Gr_G,*}\star \CF_R @>{\sim}>>
j_{\cmu,*}\underset{I}\star \CF_R @>{\on{id}_{j_{\cmu,*}}\star\phi}>> 
j_{\cmu,*}\underset{I}\star \CF^H.
\endCD
$$

\begin{lem}  \label{univ of baby Verma}
Under the above circumstances, there exists a unique map 
$\CF_{w_0}\to \CF^H$, extending $\phi$, and which intertwines
the maps \eqref{shift map} and \eqref{shift maps 2}.
\end{lem}

\ssec{}

We shall now establish the equivalence between the present definition
of $\CF_{w_0}$ and the objects defined in \cite{ABBGM}.

For a weight $\cnu\in \cLambda$ consider the inductive system of 
objects of $\on{D}(\Gr_G)_\crit\mod$, parameterized by pairs of elements 
$\cla,\cmu\in \cLambda^+\,|\, \cla-\cmu=\cnu$, and whose terms are given
by
$$j_{\cla,\Gr_G,*}\star \CF_{(V^\cmu)^*}\otimes \ell^{-\cla+\cmu}.$$

The maps in this inductive system are defined whenever two
pairs $(\cla',\cmu')$ and $(\cla,\cmu)$ are such that
$\cla'-\cla=\cmu'-\cmu=:\ceta\in \cLambda^+$, and the corresponding
map equals the composition
\begin{align*}
& j_{\cla,\Gr_G,*}\star \CF_{(V^\cmu)^*}\otimes \ell^{-\cla+\cmu}\to
j_{\cla,\Gr_G,*}\star \CF_{V^\ceta}\star \CF_{(V^\ceta)^*}\star 
\CF_{(V^\cmu)^*}\otimes \ell^{-\cla+\cmu}\to \\
&\to j_{\cla+\ceta,\Gr_G,*}\star \CF_{(V^{\cmu+\ceta})^*}\otimes
\ell^{-\cla-\ceta+(\cmu+\ceta)}.
\end{align*}

\medskip

Let $\CF'_{w_0}(\cnu)\in \on{D}(\Gr_G)_\crit\mod$ be the direct limit of 
the above system. We endow $\CF'_{w_0}:=\underset{\cnu\in \cLambda}
\oplus\, \CF'_{w_0}(\cnu)$ with the structure of an object of 
$\on{D}(\Gr_G)^\Hecke_\crit\mod$ as in Sect. 3.2.1 of \cite{ABBGM}.

\begin{prop}  \label{two approaches}
There exists a natural isomorphism 
$$\CF'_{w_0}\simeq \CF_{w_0}.$$
\end{prop}

\begin{proof}

The map $\CF_{w_0}\to \CF'_{w_0}$ is constructed using
\lemref{univ of baby Verma}, and the corresponding property
of $\CF'_{w_0}$ established in \cite{ABBGM}, Corollary 3.2.3.

To show that this map is an isomorphism, we construct
a map in the opposite direction $\CF'_{w_0}\to \CF_{w_0}$
(as mere objects of $\on{D}(\Gr_G)_\crit\mod$) as follows:

For each $\cla,\cmu\in \cLambda^+$, we let 
$j_{\cla,\Gr_G,*}\star \CF_{(V^\cmu)^*}\otimes \ell^{-\cla+\cmu}$ embed into
$j_{\cla,\Gr_G,*}\star \CF_R\otimes \ell^{-\cla}$ by means of
$$\CF_{(V^\cmu)^*}\otimes \ell^{\cmu}\hookrightarrow
\CF_{(V^\cmu)^*}\otimes \uV^\cmu\hookrightarrow \CF_R,$$ where the
second arrow is given by
$$\ell^\cmu\simeq (\uV^\cmu)^{\cN}\hookrightarrow \uV^\cmu.$$

It is straightforward to check that this gives rise to a well-defined map
from the inductive system corresponding to $\CF'_{w_0}(\cnu)$,
and that the above two maps $\CF_{w_0}\leftrightarrows \CF'_{w_0}$
are mutually inverse.

\end{proof}

\begin{cor}  \label{shift maps isom}
The maps \eqref{shift maps} 
$j_{\cmu,*}\underset{I}\star \CF_{w_0}\to \ell^\cmu\otimes \CF_{w_0}$
are isomorphisms.
\end{cor}

\begin{proof}

The assertion follows from the fact that the maps
$$j_{\cmu,*}\underset{I}\star \CF'_{w_0}(\cnu)\to \ell^{\cmu}\otimes
\CF'_{w_0}(\cnu+\cmu)$$
are easily seen to be isomorphisms.

\end{proof}

Let us now define the objects $\CF_w$ for other elements $w\in W$.
We set
$$\CF_w:=j_{w\cdot w_0,!}\underset{I}\star \CF_{w_0}.$$
In other words, if $w_0=w'\cdot w$, then
$$\CF_{w_0}\simeq j_{w',*}\underset{I}\star \CF_w.$$
{}From \propref{two approaches} it follows that $\CF_w$ are
D-modules, i.e., that no higher cohomologies appear.

\ssec{}

Let us now define the sought-after objects $\CF_w^\fZ$ of the category
$\on{D}(\Gr_G)^{\Hecke_\fZ}_\crit\mod$.

\medskip

Consider the $\cG$-torsor $\CP_{\cG,\fZ}$ over $\Spec(\fZ^\reg_\fg)$.
Recall from \secref{recol} that we have a canonical isomorphism
$\Spec(\fZ^\reg_\fg)\simeq \Op_\cg(\D)$, under which $\CP_{\cG,\fZ}$
goes over the canonical $\cG$-torsor $\CP_{\cG,\Op}$ on the
space of opers (see \cite{FG2}, Sect. 8.3, for details). Thus, we
obtain a canonical reduction of $\CP_{\cG,\fZ}$ to $\cB$ that we will
denote by $\CP_{\cB,\fZ}$.

This $\cB$-reduction defines a $\cB^-$-reduction on
$\CP_{\cG,\fZ}$. In order to define a $\cB^-$-reduction, we need to
specify for each $\cla\in \cLambda$ a line bundle, which we will
denote by $\CL^{\cla}_{w_0}$, and for each $\cla\in \cLambda^+$ a
surjective homomorphism
$$\kappa^{\cla,\fZ}:\CV^\cla_\fZ\to \CL^\cla_{w_0}.$$ These line
bundles should be equipped with isomorphisms $\CL^{\cla+\cmu}_{w_0}
\simeq \CL^\cla_{w_0} \otimes \CL^\cmu_{w_0}$, and hence give rise to
a $\cH$-torsor on $\on{Spec}(\fZ^\reg_\fg)$, which we will denote by
$\CP_{\cH,w_0}$. In addition, the maps $\kappa^{\cla,\fZ}$ should
satisfy the Pl\"ucker relations, as in \eqref{Plucker rel}. Now
observe that our $\cB$-reduction $\CP_{\cB,\fZ}$ gives rise to a
collection of compatible line subbundles $\CL^\cla$ of
$\CV^\cla_{\fZ}$. We then define $\CL^\cla_{w_0}$ as the dual of the
line bundle $\CL^{-w_0(\cla)} \hookrightarrow \CV^{-w_0(\cla)}_{\fZ}
\simeq (\CV^\cla_{\fZ})^*$.

It follows from the definition of opers (see \cite{FG2}, Sect. 1) that
the line bundle $\CL^\cla_{w_0}$ over $\Spec(\fZ^\reg_\fg)$ is
canonically isomorphic to the trivial line bundle tensored with the
one-dimensional vector space $\omega_x^{\langle
\rho,w_0(\cla)\rangle}$, where $\omega_x$ is the fiber of
$\omega_{\D}$ at the closed point $x\in \D$.

\medskip

We define the object $\wt{\CF}^\fZ_{w_0}\in
\on{D}(\Gr_G)^{\Hecke_\fZ}_\crit\mod$ as a direct sum
$$\underset{\cla\in \cLambda^+}\oplus\,
\on{Ind}^{\Hecke_\fZ}\bigl(j_{\cla,\Gr_G,*}\bigr)
\underset{\fZ^\reg_\fg}\otimes \CL_{w_0}^{-\cla}.$$ We define
$\CF^\fZ_{w_0}$ to be the quotient of $\wt{\CF}^\fZ_{w_0}$ by the same
relations as those defining $\CF_{w_0}$ as a quotient of
$\wt{\CF}_{w_0}$.

If we choose a trivialization of the $\cG$-torsor $\CP_{\cG,\fZ}$ in
such a way that $\CL^\cla_{w_0}\simeq \fZ^\reg_\fg\otimes \ell^\cla$
(such a trivialization exists), then under the equivalence
$$\on{D}(\Gr_G)^{\Hecke_\fZ}_\crit\mod\simeq
\on{D}(\Gr_G)^{\Hecke}_\crit\mod\otimes \fZ^\reg_\fg,$$ the object
$\CF^\fZ_{w_0}$ corresponds to $\CF_{w_0}$.

By construction, we have a system of maps
\begin{equation} \label{shift maps 3}
j_{\cmu,*}\underset{I}\star \CF^\fZ_{w_0}\simeq \CL_{w_0}^{\cmu}
\underset{\fZ^\reg_\fg}\otimes \CF^\fZ_{w_0},
\end{equation}
which by \corref{shift maps isom} are in fact isomorphisms.

\medskip

For other elements $w\in W$ we define
$$\CF^\fZ_w:=j_{w\cdot w_0,!}\underset{I}\star \CF^\fZ_{w_0}.$$

\ssec{}

Our present goal is to define the maps
\begin{equation} \label{sought-for maps}
\phi_w:\Gamma^{\Hecke_\fZ}(\Gr_G,\CF_w^\fZ)\to \BM_{w,\reg}\otimes 
\omega_x^{\langle 2\rho,\crho\rangle}.
\end{equation}

Since $\BM_{w,\reg}\simeq j_{w\cdot w_0,!}\underset{I}\star
\BM_{w_0,\reg}$, it is enough to define $\phi_w$ for $w=w_0$.

\medskip

Let $\CM$ be an object of $\hg_\crit\mod_\reg$. Assume that $\CM$ is
endowed with a system of maps
\begin{equation} \label{shift maps 4}
j_{\cmu,*}\underset{I}\star \CM\to \CL_{w_0}^{\cmu}
\underset{\fZ^\reg_\fg}\otimes \CM,
\end{equation}
defined for every $\cmu\in \cLambda^+$, compatible with the isomorphisms
\eqref{j's multiply} and $\CL_{w_0}^{\cmu} \underset{\fZ^\reg_\fg}\otimes
\CL_{w_0}^{\cmu'}\simeq  \CL_{w_0}^{\cmu+\cmu'}$.

Let $\phi$ be a map $\BV_\crit\to \BM$, such that for any $\cmu\in
\cLambda^+$ the diagram
\begin{equation}  \label{crucial diagram}
\CD \Gamma(\Gr_G,\CF_{V^\cmu}) @>{\beta_{V^\cmu}}>>
\CV^\cmu_\fZ\underset{\fZ^\reg_\fg} \otimes \BV_\crit
@>{\on{id}_{\CV^\cmu_\fZ}\otimes\phi}>>
\CV^\cmu_\fZ\underset{\fZ^\reg_\fg} \otimes \CM
@>{\kappa^{\cmu,\fZ}}>> \CL_{w_0}^{\cmu}
\underset{\fZ^\reg_\fg}\otimes \CM \\ @VVV & & & & @AAA \\
\Gamma(\Gr_G,j_{\cmu,\Gr_G,*}) @>{\sim}>> j_{\cmu,\Gr_G,*}\star
\BV_\crit @>{\sim}>> j_{\cmu,*}\underset{I}\star \BV_\crit
@>{\on{id}_{j_{\cmu,*}}\star \phi}>> j_{\cmu,*}\underset{I}\star \CM
\endCD
\end{equation}
is commutative.

\medskip

By the construction of $\CF_{w_0}^\fZ$, we have:

\begin{lem} \label{univ baby rep}
Under the above circumstances there exists a unique
map 
$$\Gamma^{\Hecke_\fZ}(\Gr_G,\CF_{w_0}^\fZ)\to \CM,$$
which intertwines the maps \eqref{shift maps 3} and
\eqref{shift maps 4}.
\end{lem}

Thus, to construct the map as in \eqref{sought-for maps} for $w=w_0$
we need to verify that the module
$\CM:=\BM_{w_0,\reg}\otimes \omega_x^{\langle 2\rho,\crho\rangle}$ 
possesses the required structures.

First, the map
$$\BV_\crit\to \BM_{w_0,\reg}\otimes \omega_x^{\langle
2\rho,\crho\rangle}$$ was constructed in \cite{FG2}, Sect. 7.2.

\ssec{}

To construct the data of \eqref{shift maps 4} we need to recall some
material from \cite{FG2}, Sect. 13.4. According to {\it loc. cit.}
there exists some $\cH$-torsor $\{ \cla \mapsto \CL'{}^\cla_{w_0} \}$
on $\Spec(\fZ^\reg_\fg)$ and a system of isomorphisms
$$j_{\cmu,*}\underset{I}\star \BM_{w_0,\reg}\simeq 
\CL'{}^\cla_{w_0} \underset{\fZ^\reg_\fg}\otimes \BM_{w_0,\reg}.$$

Thus, to construct the map $\phi_{w_0}$, we need to prove
the following assertion:

\begin{lem}    \label{two torsors}
There exists an isomorphism of $\cH$-torsors
$$\CL^\cmu_{w_0}\simeq \CL'{}^\cmu_{w_0}$$ which makes the diagram
\eqref{crucial diagram} commutative for $\CM:=\BM_{w_0,\reg}\otimes
\omega_x^{\langle 2\rho,\crho\rangle}$.
\end{lem}

Below we will prove this assertion by a rather explicit calculation.
In a future publication, we will discuss a more
conceptual approach. The crucial step is the following statement:

\begin{lem} \label{comp non-zero}
The composition 
$$\Gamma(\Gr_G,\CF_{V^\cmu})\to j_{\cmu,*}\underset{I}\star \BV_\crit
\overset{\on{id}_{j_{\cmu,*}}\star\phi}\to j_{\cmu,*}\underset{I}\star 
\BM_{w_0,\reg}\otimes \omega_x^{\langle 2\rho,\crho\rangle}$$
is non-zero.
\end{lem}

This proposition will be proved in \secref{proof non-zero comp}. Let
us assume it and construct the required isomorphism 
$\CL^\cmu_{w_0}\simeq \CL'{}^\cmu_{w_0}$.

\bigskip

\noindent {\em Proof of \lemref{two torsors}}.  Recall from
\cite{FG2}, Corollary 13.4.2, that there exists an isomorphism,
defined up to a scalar, $\CL^\cmu_{w_0}\simeq \CL'{}^\cmu_{w_0}$,
compatible with the action of $\Aut(\D)$. \footnote{Choosing a
coordinate $t$ on $\D$, we obtain a subgroup $\BG_m\subset \Aut(\D)$
of rescalings $t \mapsto at$.}  We will show that any choice of such
isomorphism makes the diagram \eqref{crucial diagram} commutative,
up to a non-zero scalar.

Thus, we are dealing with two non-zero maps
$$\CV^\cmu_\fZ\underset{\fZ^\reg_\fg}\otimes 
\BV_\crit \rightrightarrows \BM_{w_0,\reg}\otimes
\omega_x^{\langle \rho,w_0(\cmu)+2\crho\rangle}.$$
Recall from \cite{FG2}, Sect. 17.2, that there exists an isomorphism
$$\fZ^\reg_\fg\simeq
\Hom_{\ghat_\crit}(\BV_\crit,\BM_{w_0,\reg}\otimes \omega_x^{\langle
\rho,2\crho\rangle}),$$ compatible with the above
$\BG_m$-action. Thus, we are reduced to showing that the space
grading-preserving maps of $\fZ^\reg_\fg$-modules
$$\CV^\cmu_\fZ\to \omega_x^{\langle \rho,w_0(\cmu)\rangle}
\otimes \fZ^\reg_\fg$$
is $1$-dimensional.

However, $\CV^\cmu_\fZ$ admits a canonical filtration, whose
subquotients are isomorphic to $\omega_x^{\langle \rho,\cmu'\rangle}
\otimes \fZ^\reg_\fg$, where $\cmu'$ runs through the set weights of
$V^\cmu$ with multiplicities. For all $\cmu'\neq w_0(\cmu)$, we have
$\langle \rho,\cmu'\rangle>\langle \rho,w_0(\cmu)\rangle$. Since
the algebra $\fZ^\reg_\fg$ is non-positively graded, the above inequality
implies that the space of grading-preserving maps 
$$\omega_x^{\langle \rho,\cmu'\rangle}\otimes \fZ^\reg_\fg\to 
\omega_x^{\langle \rho,w_0(\cmu)\rangle}\otimes \fZ^\reg_\fg$$
is zero for $\cmu'\neq w_0(\cmu)$, and $1$-dimensional for
$\cmu'=w_0(\cmu)$.\qed

\ssec{Proof of \lemref{comp non-zero}}   \label{proof non-zero comp}

It is clear that if $\cmu=\cmu_1+\cmu_2$, with $\cmu_1,\cmu_2\in
\cLambda^+$, and the assertion of the proposition holds for $\cmu$,
then it also holds for $\cmu_1$. Hence it is sufficient to consider
the case of $\cmu$ that are regular.

\medskip

To prove the proposition we will use the semi-infinite cohomology
functor, denoted by $H^\semiinf\left(\fn\ppart,\fn[[t]], ?\otimes
\Psi_0\right)$, as in \cite{FG2}, Sect. 18. We will show that the
composition
\begin{align*}
&H^\semiinf\bigl(\fn\ppart,\fn[[t]], \Gamma(\Gr_G,\CF_{V^\cmu})
\otimes \Psi_0\bigr)\to H^\semiinf\bigl(\fn\ppart,\fn[[t]],
\Gamma(\Gr_G,j_{\cmu,\Gr_G,*}) \otimes \Psi_0\bigr)\to \\ &\to
H^\semiinf\bigl(\fn\ppart,\fn[[t]], \BM_{w_0,\reg}\otimes
\omega_x^{\langle 2\rho,\crho\rangle} \otimes \Psi_0\bigr)
\end{align*}
is non-zero (and, in fact, a surjection).

First, note that by \cite{FG2}, Sect. 18.3, the first arrow, i.e.,
$$H^\semiinf\bigl(\fn\ppart,\fn[[t]], \Gamma(\Gr_G,\CF_{V^\cmu})
\otimes \Psi_0\bigr)\to H^\semiinf\bigl(\fn\ppart,\fn[[t]],
\Gamma(\Gr_G,j_{\cmu,\Gr_G,*})\otimes \Psi_0\bigr)$$ is an
isomorphism. Hence, it remains to analyze the second arrow.  By
\cite{FG2}, Proposition 18.1.1, this is equivalent to analyzing the
arrow
\begin{align*}
&H^\semiinf\bigl(\fn^-\ppart,t\fn^-[[t]],
j_{w_0\cdot\crho,*}\underset{I}\star \Gamma(\Gr_G,j_{\cmu,\Gr_G,*})
\otimes \Psi_{-\crho}\bigr)\to \\
&H^\semiinf\bigl(\fn^-\ppart,t\fn^-[[t]],
j_{w_0\cdot\crho,*}\underset{I}\star \BM_{w_0,\reg}\otimes
\omega_x^{\langle 2\rho,\crho\rangle} \otimes \Psi_{-\crho}\bigr).
\end{align*}

We claim that the corresponding map
\begin{equation} \label{surj map}
j_{w_0\cdot\crho,*}\underset{I}\star
\Gamma(\Gr_G,j_{\cmu,\Gr_G,*})\simeq
j_{w_0\cdot\crho,*}\underset{I}\star j_{\cmu,*}\underset{I}\star
\BV_\crit\to j_{w_0\cdot\crho,*}\underset{I}\star
j_{\cmu,*}\underset{I}\star \BM_{w_0,\reg}\otimes \omega_x^{\langle
2\rho,\crho\rangle}
\end{equation}
is surjective for $\cmu$ regular. This would imply our claim, since the
semi-infinite cohomology functor
$H^\semiinf\bigl(\fn^-\ppart,t\fn^-[[t]],  ? \otimes \Psi_{-\crho}\bigr)$
is exact by Theorem 18.3.1 of \cite{FG2}.

\medskip

Note that $j_{w_0\cdot\crho,*}\underset{I}\star j_{\cmu,*}\simeq
j_{w_0(\cmu),*}\underset{I}\star j_{w_0\cdot\crho,*}$. Recall
from \cite{FG2}, Sect. 17.2, that we have a commutative diagram
$$ \CD j_{w_0\cdot\crho,*}\underset{I}\star \BV_\crit
@>{\on{id}_{j_{w_0\cdot\crho,*}}\star \phi}>>
j_{w_0\cdot\crho,*}\underset{I}\star \BM_{w_0,\reg}\otimes
\omega_x^{\langle 2\rho,\crho\rangle} \\ @V{\sim}VV @V{\sim}VV \\
\Gamma(\Gr_G,j_{w_0\cdot\crho,*}\underset{I}\star \delta_{1,\gr_G})
@>>> \BM_{1,\reg}\otimes \omega_x^{\langle \rho,\crho\rangle}, \endCD
$$
where the bottom arrow has the property that its cokernel, which we
denote by $\CN$, is {\it partially integrable}, i.e., it is admits
a filtration with every subquotient
integrable with respect to a sub-minimal parahoric Lie subalgebra
corresponding to some vertex $\imath$ of the Dynkin
diagram. 

Thus, the map in \eqref{surj map} can be written as
$$j_{w_0(\cmu),*}\underset{I} \star
(j_{w_0\cdot\crho,*}\underset{I}\star \BV_\crit)\to
j_{w_0(\cmu),*}\underset{I} \star (\BM_{1,\reg}\otimes
\omega_x^{\langle \rho,\crho\rangle}),$$ and since the functor
$j_{w_0(\cmu),*}\underset{I} \star ?$ is right-exact, it suffices to
show that $j_{w_0(\cmu),*}\underset{I} \star \CN$ is supported in
strictly negative cohomological degrees. In fact, we claim that this
is true for any partially integrable $I$-integrable $\hg_\crit$-module
and regular dominant coweight $\cmu$.

\medskip

Indeed, by devissage we may assume that $\CN$ is integrable with
respect to a sub-minimal parahoric corresponding to some vertex
$\imath$ of the Dynkin diagram. Then $j_{s_{\imath},*}\underset{I}
\star \CN$ lives in the cohomological degree $-1$. But since $\cmu$ is
regular, $j_{w_0(\cmu),*}\underset{I}\star j_{s_\imath,!}\simeq
j_{w_0(\cmu)\cdot s_\imath,*}$, and hence,
$$j_{w_0(\cmu),*}\underset{I}\star \CN\simeq j_{w_0(\cmu)\cdot
s_\imath,*}\underset{I}\star (j_{s_{\imath},*}\underset{I} \star
\CN),$$ and our assertion follows from the fact that the functor of
convolution with $j_{w_0(\cmu)\cdot s_\imath,*}$ is right-exact.\qed

\ssec{Proof of \corref{F non-zero} and completion of the
proof of \thmref{main}}

Thus, we have proved \lemref{comp non-zero} and therefore \lemref{two
torsors}. By \lemref{univ baby rep}, this implies that we have a
canonical map
$$
\phi_{w_0}: \Gamma^{\Hecke_\fZ}(\Gr_G,\CF_{w_0}^\fZ)\to
\BM_{w,\reg}\otimes \omega_x^{\langle 2\rho,\crho\rangle}.
$$
According to the remark after formula \eqref{sought-for maps}, we then
obtain maps
$$
\phi_w:\Gamma^{\Hecke_\fZ}(\Gr_G,\CF_w^\fZ)\to \BM_{w,\reg}\otimes
\omega_x^{\langle 2\rho,\crho\rangle}
$$
for all $w \in W$ (as in formula \eqref{sought-for maps}).

\begin{prop}  \label{one of the maps surj}
The map
$$\phi_1:\Gamma^{\Hecke_\fZ}(\Gr_G,\CF_1^\fZ)\to \BM_{1,\reg}\otimes 
\omega_x^{\langle 2\rho,\crho\rangle}$$
is surjective.
\end{prop}

Since the functors $j_{w,*}$ are right-exact, this proposition implies
that the same surjectivity assertion holds for all $w\in W$. Hence,
\propref{one of the maps surj} implies \corref{F non-zero} and
\thmref{main}.

\bigskip

\noindent {\em Proof of \propref{one of the maps surj}}. For $\cla$,
such that $\cla-\crho$ is dominant and regular, let us consider the
map
$$j_{w_0\cdot \cla,*}\underset{I}\star \CF_{R_\fZ}\underset{\fZ^\reg_\fg}
\otimes \CL_{w_0}^{-\cla}\simeq j_{w_0,!}\underset{I}\star
j_{\cla,*}\underset{I}\star \CF_{R_\fZ}\underset{\fZ^\reg_\fg} \otimes
\CL_{w_0}^{-\cla}\to j_{w_0,!}\underset{I}\star \wt{\CF}^\fZ_{w_0}\to
j_{w_0,!}\underset{I}\star \CF^\fZ_{w_0}\simeq \CF^\fZ_1,$$ and the
resulting map
$$j_{w_0\cdot \cla,*}\underset{I}\star \BV_\crit
\underset{\fZ^\reg_\fg} \otimes \CL_{w_0}^{-\cla}\to
\Gamma^{\Hecke_\fZ}(\Gr_G,\CF_1^\fZ)\overset{\phi_1}\to
\BM_{1,\reg}\otimes \omega_x^{\langle 2\rho,\crho\rangle}.$$

By construction, this map is obtained by applying the functor
$j_{w_0\cdot \cla,*}\underset{I}\star ?$ to the map
$$\BV_\crit\to \BM_{w_0,\reg}\otimes \omega_x^{\langle
2\rho,\crho\rangle},$$ and it coincides with the map from \eqref{surj
map} for $\cmu=\cla-\crho$.  Hence, it is surjective by \secref{proof
non-zero comp}.\qed

\ssec{Completion of the proof of \thmref{get Wakimoto}}

Thus, the proof of \thmref{main} is complete. Let us now finish the
proof of the fact that the morphisms $\phi_w$ are actually
isomorphisms and hence complete our proof of \thmref{get
Wakimoto}. Clearly, it is enough to do so for just one element of
$W$. We shall give two proofs.

\medskip

\noindent{\bf Proof 1.}  This argument will rely on \thmref{main}. We
will analyze the map $\phi_{w_0}$.  By \cite{ABBGM}, Proposition
3.2.5, the canonical map $\CF_R \to \CF_{w_0}$ identifies
$\on{Ind}^\Hecke(\delta_{1,\Gr_G})$ with the co-socle of
$\CF_{w_0}$. Hence $\Gamma^{\Hecke_\fZ}(\Gr_G,\CF^\fZ_{w_0})$ does not
have sub-objects whose intersection with $\BV_\crit =
\Gamma^{\Hecke_\fZ}(\Gr_G,R_{\fZ})$ is zero.

Therefore, to prove the injectivity of the map $\phi_{w_0}$, it is
enough to show that the composition
$$\BV_\crit\simeq\Gamma^{\Hecke_\fZ}(\Gr_G,R_\fZ)\to
\Gamma^{\Hecke_\fZ}(\Gr_G,\CF^\fZ_{w_0})\overset{\phi_{w_0}}\to
\BM_{w_0,\reg} \otimes \omega_x^{\langle 2\rho,\crho\rangle}$$ is
injective. However, the latter map is, by construction, the map
$\BV_\crit\to \BM_{w_0,\reg}\otimes \omega_x^{\langle
2\rho,\crho\rangle}$ of \cite{FG2}, Sect. 17.2, which was injective by
definition.

\medskip

\noindent{\bf Proof 2.} This argument will be independent of
\thmref{main},(2). We will analyze the map $\phi_1$. We have a
canonical map
$$\IC_{w_0\cdot \crho,\Gr}\star \CF_{R_\fZ}\underset{\fZ^\reg_\fg}
\otimes \CL_{w_0}^{-\crho}\to
j_{w_0\cdot \crho}\underset{I}\star \CF_{R_\fZ}\underset{\fZ^\reg_\fg}
\otimes \CL_{w_0}^{-\crho}\to \CF^\fZ_1,$$
and by \cite{ABBGM}, Propositions 3.2.6 and 3.2.10, its cokernel
is partially integrable.

The composition 
\begin{align*}
&\Gamma(\Gr_G, \IC_{w_0\cdot \crho,\Gr})
\underset{\fZ^\reg_\fg}\otimes \CL_{w_0}^{-\crho}\simeq
\Gamma^{\Hecke_\fZ}\bigl(\Gr_G,\IC_{w_0\cdot \crho,\Gr}
\star \CF_{R_\fZ}\underset{\fZ^\reg_\fg}\otimes \CL_{w_0}^{-\crho}\bigr)
\simeq \Gamma^{\Hecke_\fZ}\bigl(\Gr_G,\CF^\fZ_1)
\overset{\phi_1}\to \\
&\to \BM_{1,\reg}\underset{\fZ^\reg_\fg}\otimes 
\omega_x^{\langle 2\rho,\crho\rangle}
\end{align*}
comes from the map
$$\Gamma(\Gr_G, \IC_{w_0\cdot \crho,\Gr})\to
\BM_{1,\reg}\underset{\fZ^\reg_\fg}\otimes \omega_x^{\langle
\rho,\crho\rangle},$$ of \cite{FG2}, Sect. 17.3, which is injective by
{\it loc.cit.}

Hence, the kernel of the map $\phi_1$ is partially integrable. But we claim
that $\Gamma^{\Hecke_\fZ}\bigl(\Gr_G,\CF^\fZ_1)$ admits no partially
integrable submodules.

Indeed, suppose that $\CN$ is a submodule of
$\Gamma^{\Hecke_\fZ}\bigl(\Gr_G,\CF^\fZ_1)$, integrable with respect
to a sub-minimal parahoric, corresponding to a vertex $\imath$ of the
Dynkin diagram. Since the functor $j_{s_\imath,*}\underset{I}\star$ is
invertible on the derived category, we would obtain a non-zero map:
$$j_{s_\imath,*}\underset{I}\star \CN\to \on{L}
\Gamma^{\Hecke_\fZ}\bigl(\Gr_G,j_{s_\imath,*}
\underset{I}\star\CF^\fZ_1).$$

But the LHS is supported in the cohomological degrees $<0$, and the
RHS is acyclic away from cohomological degree $0$. \footnote{Here we
are relying on part (1) of \thmref{main}, which was proved
independently.} This is a contradiction.

This completes the proof of \thmref{get Wakimoto}.

\section{Appendix: an equivalence at the negative level}

\ssec{}

Let $\kappa$ be a negative level, i.e., $\kappa=k\cdot \kappa_{\can}$
with $k+h^\vee\notin \BQ^{\geq 0}$.

Let $\wt\Fl_G$ be the enhanced affine flag scheme, i.e, $G\ppart/I^0$, and
let $\on{D}(\wt\Fl_G)_\kappa\mod$ be the corresponding
category of twisted D-modules. 

Note that $\wt\Fl_G$ is acted on by the group $I/I^0\simeq H$ by right
multiplication. Let us denote by $\on{D}(\wt\Fl_G)_\kappa\mod^{H,w}$
the corresponding category of weakly H-equivariant objects of
$\on{D}(\wt\Fl_G)_\kappa\mod$ (see \cite{FG2}, Sect. 20.2). 

For an object $\CF\in \on{D}(\wt\Fl_G)_\kappa\mod^{H,w}$, consider
$\Gamma(\wt\Fl_G,\CF)\in \hg_\kappa\mod$.
The weak $H$-equivariant structure on 
$\CF$ endows $\Gamma(\wt\Fl_G,\CF)$ with a commuting action of $H$.
We let 
$$\Gamma^H:\on{D}(\wt\Fl_G)_\kappa\mod^{H,w}\to \hg_\kappa\mod$$
to be the composition of $\Gamma(\wt\Fl_G,\cdot)$, followed by the
functor of $H$-invariants. 

\medskip

Recall from \cite{FG2}, Sect. 20.4, that 
every object of $\on{D}(\wt\Fl_G)_\kappa\mod^{H,w}$ carries
a canonical action of $\Sym(\fh)$ by endomorphism, denoted
$a^\sharp$.

For $\lambda\in \fh^*$ let 
$$\on{D}(\wt\Fl_G)_\kappa\mod^{H,\lambda}\subset
\on{D}(\wt\Fl_G)_\kappa\mod^{H,w,\lambda}$$ be the full subcategories
of $\on{D}(\wt\Fl_G)_\kappa\mod^{H,w}$, corresponding to the condition
that $a^\sharp(h)=\lambda(h)$ for $h\in \fh$ in the former case, and
that $a^\sharp(h)-\lambda(h)$ acts locally nilpotently in the latter.
Since the group $H$ is connected, both of these categories are full
subcategories in $\on{D}(\wt\Fl_G)_\kappa\mod$.

\medskip

We let $D(\on{D}(\wt\Fl_G)_\kappa\mod)^{H,w,\lambda}\subset
D(\on{D}(\wt\Fl_G)_\kappa\mod)$ be the full subcategory consisting of
complexes, whose cohomologies belong to
$\on{D}(\wt\Fl_G)_\kappa\mod^{H,w,\lambda}$.
It is easy to see that the functor $\Gamma^H$, restricted to
$\on{D}(\wt\Fl_G)_\kappa\mod^{H,w,\lambda}$, extends to a functor
$$\on{R}\Gamma^H:D^+(\on{D}(\wt\Fl_G)_\kappa\mod)^{H,w,\lambda}\to
D^+(\hg_\kappa\mod).$$

\medskip

Assume now that $\lambda$ satisfies the following conditions:
$$
\begin{cases}
& \langle \lambda+\rho ,\check\alpha\rangle\notin \BZ^{\geq 0} \text{
for } \alpha \in \Delta^+ \\ &\pm \langle \lambda+\rho
,\check\alpha\rangle +2n\frac{k+h^\vee}{\kappa_{\can}(\alpha,\alpha)}
\notin \BZ^{\geq 0} \text{ for } \alpha \in \Delta^+ \text{ and } n\in
\BZ^{>0}.
\end{cases}
$$

Following \cite{BD}, Sect. 7.15, we will prove:

\begin{thm} \label{weak local, neg}   \hfill

\smallskip

\noindent{\em (1)} For $\CF\in
\on{D}(\wt\Fl_G)_\kappa\mod^{H,w,\lambda}$ the higher cohomologies
$\on{R}^i\Gamma^H(\wt\Fl_G,\CF)$, $i>0$, vanish.

\smallskip

\noindent{\em (2)}
The resulting functor $\on{R}\Gamma^H:
D^b(\on{D}(\wt\Fl_G)_\kappa\mod)^{H,w,\lambda}\to
D^b(\hg_\kappa\mod)$
is fully-faithful.

\end{thm}

\ssec{}

Let $\on{D}(\wt\Fl_G)_\kappa\mod^{I^0,H,w,\lambda}\subset
\on{D}(\wt\Fl_G)_\kappa\mod)^{H,w,\lambda}$ be the full subcategory,
consisting of twisted D-modules, equivariant with respect to the
$I^0$-action on the left. Our present goal is to describe its image
under the above functor $\Gamma$.

Consider the category $\CO_{\aff}:=\hg_\kappa\mod^{I^0}$. This is a
version of the category $\CO$ for the affine Lie algebra
$\ghat_\kappa$. Its standard (resp., co-standard, irreducible) objects
are numbered by weights $\mu\in \fh^*$, and will be denoted by
$M_{\kappa,\mu}$ (resp., $M^\vee_{\kappa,\mu}$, $L_{\kappa,\mu}$).
Since $\kappa$ was assumed to be negative, every finitely generated
object of $\CO_{\aff}$ has finite length.

The extended affine Weyl group $W_{\aff}:=W\ltimes \Lambda$ acts on
$\fh^*$, with $w\in W\subset W_{\aff}$ acting as
$$w\cdot \mu=w(\mu+\rho)-\rho,$$ and $\cla\in \cLambda\subset
W_{\aff}$ by the translation by means of
$(\kappa-\kappa_{\crit})(\cla,\cdot)\in \fh^*$.

For a $W_{\aff}$-orbit $\upsilon$ in $\fh^*$ let
$(\CO_{\aff})_{\upsilon}$ be the full-subcategory of $\CO_{\aff}$,
consisting objects that admit a filtration, such that all subquotients
are isomorphic to $L_{\kappa,\lambda}$ with $\lambda\in \upsilon$.

The following assertion is known as the linkage principle (see
\cite{DGK}):

\begin{prop}
The category $\CO_{\aff}$ is the direct sum over the orbits $\upsilon$
of the subcategories $(\CO_{\aff})_{\upsilon}$.
\end{prop}

For $\lambda$ as in \thmref{weak local, neg} let $\upsilon(\lambda)$
be the $W_\aff$-orbit of $\lambda$. (Note that by assumption, the
stabilizer of $\lambda$ in $W_\aff$ is trivial.)

We shall prove the following: \footnote{This theorem is not due to the
authors of the present paper. The proof that we present is a
combination of arguments from \cite{BD}, Sect. 7.15, and \cite{KT}.}

\medskip

\begin{thm}   \label{KTthm}
The functor $\Gamma^H$ defines an equivalence
$$\on{D}(\wt\Fl_G)_\kappa\mod)^{I^0,H,w,\lambda}\to
(\CO_{\aff})_{\upsilon(\lambda)}.$$
\end{thm}

\ssec{Proofs}

To prove point (1) of \thmref{weak local, neg}, it suffices to show
that $\on{R}^i\Gamma^H(\wt\Fl_G,\CF)=0$ for $\CF\in
\on{D}(\wt\Fl_G)_{\kappa}\mod^{H,\lambda}$ and $i>0$. However, this
follows immediately from \cite{BD}, Theorem 15.7.6.

\medskip

To prove point (2) of \thmref{weak local, neg} and \thmref{KTthm} we
shall rely on the following explicit computation, performed in
\cite{KT}:

For an element $\wt{w}\in W_{\aff}$ let $j_{\wt{w},*,\lambda}\in
\on{D}(\wt\Fl_G)_{\kappa}\mod^{I^0,H,\lambda}$ (resp.,
$j_{\wt{w},!,\lambda}$) be the *-extension (resp, !-extension) of the
unique $I^0$-equivariant irreducible twisted D-module on the preimage
of the corresponding $I^0$-orbit in $\Fl_G$. We have:

\begin{thm} \label{dual Verma} We have:
$$\Gamma(\Fl_G,j_{\wt{w},*,\lambda})\simeq M^\vee_{\kappa, \wt{w}\cdot
0} \text{ and } \Gamma(\Fl_G,j_{\wt{w},!,\lambda})\simeq M_{\kappa,
\wt{w}\cdot 0}.$$
\end{thm}

Let us now proceed with the proof of \thmref{weak local,
neg}(2). Clearly, it is enough to show that for two finitely generated
objects $\CF,\CF_1\in \on{D}(\wt\Fl_G)_{\kappa}\mod^{H,\lambda}$ the
map
$$\on{R}\on{Hom}_{D(\on{D}(\wt\Fl_G)_\kappa\mod)^{H,\lambda}}(\CF,
\CF_1) \to \on{R} \on{Hom}_{D(\hg_\crit\mod)}(\Gamma^H(\wt\Fl_G,
\CF),\Gamma^H(\wt\Fl_G,\CF))$$ is an isomorphism.

By adjunction (see \cite{FG2}, Sect. 22.1), the latter is equivalent
to the map
\begin{align*}
&\on{R}\on{Hom}_{D(\on{D}(\wt\Fl_G)_\kappa\mod)^{I,\lambda}}
(j_{1,!,\lambda},\CF^{op}\star \CF_1)\to \\
&\to\on{R}\on{Hom}_{D(\hg_\crit\mod)^{I,\lambda}}(
\Gamma^H(\wt\Fl_G,j_{1,!,\lambda}), \on{R}\Gamma^H(\Fl_G,\CF^{op}
\star \CF_1))
\end{align*}
being an isomorphism, where $\CF^{op}\in
\on{D}(G\ppart/K)\mod^{I,\lambda}$ is the dual D-module, where $K$ is
a sufficiently small open-compact subgroup of $G[[t]]$.

Using the stratification of $\wt\Fl_G$ by $I$-orbits, we can replace
$\CF^{op}\star \CF_1$ by its Cousin complex. In other words, it is
sufficient to show that
$$\on{R}\on{Hom}_{D(\on{D}(\wt\Fl_G)_\kappa\mod)^I}(j_{1,!,\lambda},
j_{\wt{w},\lambda,*})\to
\on{R}\on{Hom}_{D(\hg_\crit\mod)^I}(\Gamma^H(\wt\Fl_G,j_{1,!,\lambda}),
\Gamma^H(\wt\Fl_G,j_{\wt{w},\lambda,*}))$$ is an isomorphism, for all
$\wt{w}$ such that $j_{\wt{w},\lambda,*}$ is
$(I,\lambda)$-equivariant.

Note that the LHS is $0$ unless $\wt{w}=0$, and is isomorphic to $\BC$
in the latter case. Hence, taking into account \thmref{dual Verma}, it
remains to prove the following:

\begin{lem} \hfill

\smallskip

\noindent{\em (1)}
$\on{R}\on{Hom}_{D(\hg_\crit\mod)^{I,\lambda}}(M_{\kappa,\lambda},
M^\vee_{\kappa,\mu})=0$ for $\lambda\neq \mu\in \fh^*$ but such that
$M^\vee_{\kappa,\mu}\in \hg_\crit\mod^{I,\lambda}$ is
$(I,\lambda)$-equivariant.

\smallskip

\noindent{\em (2)} The map $\BC\to
\on{R}\on{Hom}_{D(\hg_\crit\mod)^{I,\lambda}}(M_{\kappa,\lambda},
M^\vee_{\kappa,\lambda})$ is an isomorphism.

\end{lem}

\begin{proof}

For any $\CM\in \hg_\crit\mod^{I,\lambda}$,
$$\on{R}\on{Hom}_{D(\hg_\crit\mod)^{I,\lambda}}(M_{\kappa,\lambda},
\CM)\simeq \on{R}\on{Hom}_{I\mod}(\BC,\CM\otimes \BC^{-\lambda}).$$

Since $M^\vee_{\kappa,\mu}$ is co-free with respect to $I^0$,
we obtain
$$\on{R}\on{Hom}_{I\mod}(\BC,M^\vee_{\kappa,\mu}\otimes
\BC^{-\lambda})\simeq \on{R}\on{Hom}_{H\mod}(\BC,\BC^\mu\otimes
\BC^{-\lambda}),$$ implying the first assertion of the lemma.

Similarly,
$$\on{R}\on{Hom}_{D(\hg_\crit\mod)^I}(M_{\kappa,\lambda},M^\vee_{\kappa,\lambda})\simeq
\on{R}\on{Hom}_{I\mod}(\BC,M^\vee_{\kappa,\lambda})\simeq 
\on{R}\on{Hom}_{H\mod}(\BC,\BC)\simeq \BC,$$
implying the second assertion.

\end{proof}

Finally, let us prove \thmref{KTthm}. Taking into account \thmref{weak
local, neg}, and using Lemmas \ref{adjoint} and \ref{kappa adjoint},
it remains to show that for every $\CM\in
(\CO_{\aff})_{\upsilon(\lambda)}$ there exists an object $\CF\in
\on{D}(\wt\Fl_G)\mod^{I^0,H,w,\lambda}$ with non-zero map
$$\Gamma^H(\wt\Fl_G,\CF)\to \CM.$$

It is clear that for every $\CM\in (\CO_{\aff})_{\upsilon(\lambda)}$
there exists a Verma module $M_{\kappa,\mu}\in
(\CO_{\aff})_{\upsilon(\lambda)}$ with a non-zero map
$M_{\kappa,\mu}\to \CM$. Hence, the required property follows from
\thmref{dual Verma}.

\end{document}